\documentclass[12pt]{amsart}

\usepackage{verbatim}
\usepackage{cite}
\usepackage{adjustbox}
\usepackage{amssymb}
\usepackage{amsmath}
\usepackage{amsthm}
\usepackage{amsfonts}
\usepackage{float}
\usepackage{graphicx}%
\graphicspath{{./}}
\usepackage{cancel}
\usepackage{comment}
\usepackage{bbm}
\usepackage{xcolor}
\usepackage[toc,page]{appendix}
\usepackage{mathtools}
\usepackage{soul}

\usepackage[plain]{fullpage}
\usepackage[stretch=10,nopatch=footnote]{microtype}
\usepackage{todonotes}
\usepackage{bm}
\usepackage{enumitem}
\usepackage[hyperfootnotes=false]{hyperref}

\newcommand{\samuel}[1]{{\color{teal} Samuel says: #1}}

\numberwithin{equation}{section}
\theoremstyle{plain}
\newtheorem{theorem}{Theorem}[section]
\newtheorem{proposition}[theorem]{Proposition}
\newtheorem{claim}[theorem]{Claim}
\newtheorem{lemma}[theorem]{Lemma}
\newtheorem{corollary}[theorem]{Corollary}

\theoremstyle{definition}
\newtheorem{definition}[theorem]{Definition}
\newtheorem*{defs*}{Definition}
\newtheorem{conj}[theorem]{Conjecture}

\theoremstyle{remark}
\newtheorem{remark}[theorem]{Remark}
\newtheorem*{ass*}{Assumption}
\newtheorem*{ack*}{Acknowledgements}

\usepackage{setspace} 

\title{Geodesics on Grushin spaces}

\begin{document}
\author[Albert]{Michael Albert{$^{\dag }$}}
\address{Department of Mathematics\\
    University of Connecticut\\
    Storrs, CT 06269,  U.S.A.}
\email{michael.2.albert@uconn.edu}

\author[Borza]{Samu\"{e}l Borza{$^{\ddag}$}}
\thanks{{$^{\ddag}$}
    This research was funded in part by the Austrian Science Fund (FWF) [10.55776/EFP6].}
\address{Faculty of Mathematics, University of Vienna, Oskar-Morgenstern-Platz 1,
    1090,
    Austria}
\email{samuel.borza@univie.ac.at}

\author[Gordina]{Maria Gordina{$^{\dag }$}}
\thanks{{$^{\dag}$} Research was supported in part by NSF grants DMS-2246549 and DMS-2154047 while both MA and MG were at the University of Connecticut}
\address{ Department of Mathematics\\
    University of Rochester\\
    Rochester, NY 14627,  U.S.A.}
\email{mgordina@ur.rochester.edu}



\date{\today \ \emph{File:\jobname{.tex}}}

\begin{abstract}
    We consider higher-dimensional generalizations of the $\alpha$-Grushin plane, focusing on the problem of classification of geodesics that minimize length, also known as \textit{optimal synthesis}. Solving Hamilton's equations on these spaces using the calculus of generalized trigonometric functions, we obtain explicit conjugate times for geodesics starting at a Riemannian point. We propose a conjectured cut time $\tau=\min\{\tau_j\}$ obtained as the minimum of several candidates, each deriving from the symmetries of components of the Hamiltonian flow. We prove that it provides a lower bound on the first conjugate time, a key step in the extended Hadamard technique. In the three-dimensional case, we combine this method with a density argument to establish the conjecture and obtain the full optimal synthesis.
\end{abstract}

\maketitle

\tableofcontents

\section{Introduction}
Since the seminal work of Grushin and Baouendi on hypoelliptic operators \cite{Grushin1970,Baouendi1967}, the geometric structures nowadays commonly refered to as \emph{Grushin spaces} have been widely studied. The most well-known is the Grushin plane, the sub-Riemannian structure on $\mathbb{R}^2$ generated by the global family of vector fields
\[
    X = \partial_x, \qquad Y = x \partial_y.
\]
A variant, the \emph{$\alpha$-Grushin plane}, is a obtained by replacing the vector field $Y$ with $Y_\alpha = |x|^\alpha \partial_y$ for $\alpha\geq 0$.

In this paper we consider a higher dimensional generalization of the $\alpha$-Grushin plane. Given $\alpha = (\alpha_1, \dots, \alpha_n) \in \mathbb{N}^{n}$, we call \emph{$\alpha$-Grushin space}, the space obtained by equipping $\mathbb{R}^{n + 1}$ with the sub-Riemannian structure generated by the global family of vector fields
\[
    X_1 = \partial_{x_1} \ , \quad X_2 = x_1^{\alpha_1} \partial_{x_2} \ , \quad X_3 = x_1^{\alpha_1} x_2^{\alpha_2} \partial_{x_3} \ , \quad \dots \ , \qquad X_{n + 1} = x_1^{\alpha_1} x_2^{\alpha_2} \cdots x_n^{\alpha_n} \partial_{x_{n+1}}.
\]
Our work initiates the study of the \emph{optimal synthesis} of this space, that is to say, the explicit construction of the complete set of geodesics starting from any point, as well as their associated \emph{cut} and \emph{conjugate times}.

The higher-dimensional $\alpha$-Grushin space that we consider has been studied by a number of authors including \cite{FerrariFranchi2000}, and in the context of $p$-Laplacians, viscosity and other problems in PDEs such as \cite{Bieske2005LipschitzEO, BieskeGong2006}. A geometric perspective on this model of Grushin spaces has been presented in \cite{WuJang-Mei2015}, while \cite{TurcanuUdriste2017} considered this model from the perspective of stochastic accessibility. An important application of Grushin operators and Grushin geodesics was found in the context of metamaterial arrays in \cite{GreenleafKettunenKurylevLassasUhlmann2018}. We refer to \cite{BarilariBoscainNeel2012, BoscainPrandi2016, BoscainNeel2020, BorzaPhDThesis2021, WuJang-Mei2015, KogojLanconelli2012} for studies of spaces equipped with geometries induced by Grushin-type vector fields and the associated Grushin operators in the non-smooth, non-H\"{o}rmander setting.

A geometric feature of the classical Grushin plane is that its Riemannian region has
negative Gaussian curvature, with curvature blowing up near the singular set, while
conjugate points nevertheless occur after geodesics cross the singular line. This
phenomenon, discussed in the almost-Riemannian literature, shows that the usual
Riemannian intuition that negative curvature prevents conjugate points is not reliable
in the presence of a sub-Riemannian distribution; see, for instance,
\cite{AgrachevBoscainSigalotti2008,BoscainLaurent2013}. This phenomenon persists in the higher dimensional Grushin spaces that we consider, and conjugacy past the singular set crossing will be detected despite strictly negative sectional, Ricci and scalar curvatures on the Riemannian portion of the space. See Theorem \ref{Con Times 3D}.

The three-dimensional case is related, but subtly different from the local theory of almost-Riemannian
structures developed by Boscain--Charlot--Gaye--Mason
\cite{BoscainCharlotGayeMason2015}. Indeed we remark explicitly that this work is not directly applicable to 3D Grushin spaces. Their generic
hypotheses lead to singular points of rank at least \(2\) and both their ``Type 1" and ``Type 2" inherit bracket assumptions that are not satisfied in our setting. At certain singular points in Grushin spaces, the missing directions are
recovered only after potentially higher-order brackets.

We now place the present work in the context of known optimal-synthesis results. The
Heisenberg group was the first fundamental model where the sub-Riemannian wave front,
sphere, cut locus, and conjugate locus were understood explicitly, but this understanding
did not come from a single source. Gaveau's heat-kernel work \cite{Gaveau1977a},
Brockett's control-theoretic analysis \cite{Brockett1982}, Nachman's study of the wave
equation \cite{Nachman1982}, and the work of Vershik--Gershkovich on nonholonomic
variational problems \cite{VershikGershkovich1988} all contributed to the development of
this model. The optimal synthesis has since been completed for several other step-two
structures: step-two contact nilpotent structures \cite{AgrachevBarilariBoscain2012},
step-two corank-two structures \cite{BarilariBoscainGauthier2012}, and, more recently,
Reiter--Heisenberg groups \cite{MontanariMorbidelli2024}. The free step-two Carnot
groups remain more subtle: an early conjectural description of the cut locus
\cite{Myasnichenko2002, Myasnichenko2006, Monroy-PerezAnzaldo-Meneses2006} was later
disproved in \cite{RizziSerres2017}, and the full problem remains open.

Left-invariant sub-Riemannian structures on three-dimensional Lie groups were classified
in \cite{AgrachevBarilari2012}, and the cut locus has been studied in several of these
models, including \(\operatorname{SU}(2)\), \(\operatorname{SO}(3)\), and
\(\operatorname{SL}(2)\) \cite{BoscainRossi2008a}; the group
\(\operatorname{SE}(2)\) of motions of the plane
\cite{MoiseevSachkov2010, Sachkov2010a, Sachkov2011}; and the group
\(\operatorname{SH}(2)\) of hyperbolic motions of the plane
\cite{ButtSachkovBhatti2014, ButtSachkovBhatti2016, ButtSachkovBhatti2017}. In higher
step, explicit optimal synthesis is substantially more difficult, but it is known in
important examples such as the Martinet flat structure
\cite{AgrachevBonnardChybaKupka1997}, the Engel group
\cite{ArdentovSachkov2011, ArdentovSachkov2013, ArdentovSachkov2015,
    ArdentovSachkov2017}, and the Cartan group
\cite{Sachkov2021a, ArdentovHakavuori2022}.

For Grushin-type geometries, the analytic origins go back at least to Baouendi's work on
degenerate elliptic operators \cite{Baouendi1967} and to Grushin's hypoelliptic
operators \cite{Grushin1970}. Bellaïche's account of sub-Riemannian geometry contains
early pictures and discussion of Grushin geodesics and spheres \cite{Bellaiche1996}.
For the classical Grushin plane, the cut locus is described in full detail using the Extended Hadamard technique in
\cite[Chapter~13]{AgrachevBarilariBoscainBook2020}, while the general
\(\alpha\)-Grushin plane was treated in \cite{Borza2022}; geodesics in related
Grushin-type spaces were also studied in \cite{ChangLi2012,changli2014}. From the
almost-Riemannian viewpoint, Agrachev--Boscain--Sigalotti studied the conjugate locus
from a nonsingular point and emphasized the relation between singularities, curvature,
and conjugate points \cite{AgrachevBoscainSigalotti2008}. The associated
Laplace--Beltrami operator, and the role of the singular set for diffusion and quantum
confinement, were studied by Boscain--Laurent \cite{BoscainLaurent2013}.


The study of the optimal synthesis is relevant to many applications. For example, it is important for generalizations of curvature on sub-Riemannian manifolds, where the absence of a Levi-Civita connection makes defining curvature difficult. Various notions of curvature in sub-Riemannian geometry have appeared over the last several years. The approach of Sturm, Lott, Vilani and Ohta introduces \emph{synthetic} curvature-dimension bounds with  optimal transport theory, the goal being to generalize the Riemannian condition ``\textit{$\operatorname{Ric} \geqslant K$ and $n \leqslant N$}'' to metric measure spaces. The \emph{measure contraction property} $\mathsf{MCP}(K, N)$, for instance, has been shown to hold for many sub-Riemannian manifolds. With an explicit expression for the cut times, one can optimize parameters $K$ and $N$ in the $\mathsf{MCP}(K, N)$ and obtain the sharpest information that this synthetic notion of curvature can provide. This was done for the Heisenberg group, the $\alpha$-Grushin plane, and $H$-type groups \cite{Juillet2009,BarilariRizzi2018,Borza2022,Rizzi2016}, among others. It is also the optimal synthesis of the Martinet flat structure \cite{AgrachevBonnardChybaKupka1997} that allowed the authors of \cite{borzaMCPfailure} to show that the $\mathsf{MCP}(K, N)$ actually fails in large families of sub-Riemannian manifolds.
Optimal synthesis can also be applied to \emph{geodesic random walks} on sub-Riemannian manifolds, as considered in \cite{GordinaLaetsch2017, AgrachevBoscainNeelRizzi2018, BoscainNeelRizzi2017, GordinaMelcherMikulincerWang2024}.

The paper is organized as follows. In Section~2 we recall the sub-Riemannian
preliminaries used throughout the paper. In Section~3 we define $\alpha$-Grushin spaces and derive their Hamiltonian system. Section~4 introduces the generalized
trigonometric functions needed to solve Hamilton's equations explicitly. In
Section~5 we classify abnormal minimizers and show that, although higher-dimensional $\alpha$-Grushin spaces are not ideal in general, their abnormal minimizers are not
strictly abnormal. Section~6 gives the determinant factorization for the exponential
map, proves a lower bound for the first conjugate time, and formulates the main cut-time
conjecture. In Section~7 we prove the conjecture for \(\mathbb G^3_{(\alpha,\beta)}\)
at Riemannian points using the extended Hadamard technique. Finally, Section~8 extends
the optimal synthesis in $\mathbb{G}^3_{(\alpha,\beta)}$ to singular initial points by density.

We make use of the numerical ODE solver capabilities of \emph{Mathematica} (Wolfram Research, Champaign, IL, USA) to generate simulations of geodesics and the cut loci for various points in $\mathbb{G}^3_{(\alpha,\beta)}$\footnote{\textbf{Data Availability Statement.}
    The code used in this paper is available online in a Github repository: \url{https://github.com/StevieMike/Geodesics-on-Grushin-Spaces-Mathematica-Notebook}}.

\begin{ack*} The authors would like to thank Matthew Badger for his invaluable guidance on this project.
\end{ack*}

\section{Preliminaries on sub-Riemannian geometry}\label{s.subRiemBasics}

For a comprehensive introduction to the basic theory of sub-Riemannian manifolds, see the textbook \cite{AgrachevBarilariBoscainBook2020}. We will work in the setting of a globally smooth, \emph{bracket generating} family of vector fields $\mathcal{F}=\{Z_1,...,Z_m\}$ on a smooth manifold $M$ of dimension $N$. For our purposes a \emph{sub-Riemannian structure} is a pair $(M,\mathcal{F})$.

We define a family of subspaces, called the \emph{admissible distribution} \( \mathcal{D} = \{\mathcal{D}_q\}_{q\in M} \), where $\mathcal{D}_q=\operatorname{span}\{Z_j(q):1\leq j\leq m\}$. Contrary to certain treatments, we do not impose that the frame $\{Z_1(q),...,Z_m(q)\}$ has constant rank, and therefore $\mathcal{D}=\{\mathcal{D}_q\}_{q\in M}\subset TM$ may not have the structure of a sub-bundle.

A Lipschitz (in charts) curve $\gamma : [0, T] \to M$ is said to be \emph{admissible} (or \emph{horizontal}) if there exists $u \in L^\infty([0,T],\mathbb{R}^m)$, called a \emph{control}, such that
\begin{align}\label{Admissible}
    \gamma^{\prime}(t)=\sum_{j=1}^m u_j(t)Z_j(\gamma(t)), \ \qquad  \text{for a.e. } t\in [ 0, T].
\end{align}
For fixed $t\in[0,T]$ such that \eqref{Admissible} holds, let $u^{\ast}(t)$ be the unique minimizer of $\lvert v \rvert$ among all vectors $v\in \mathbb{R}^m$ satisfying $\gamma'(t)=\sum_{j=1}^mv_jZ_j(\gamma(t))$. We refer to \cite[Chapter~3]{AgrachevBarilariBoscainBook2020} for a proof that the resulting \emph{minimal control} $t\mapsto u^{\ast}(t)$ is measurable and essentially bounded. The \emph{(sub-Riemannian) length} of an admissible curve is defined as
\begin{align*}
    \ell(\gamma):=\int_0^T\lvert u^{\ast}(t)\rvert \text{d}t.
\end{align*}
Any admissible curve admits a Lipschitz reparametrization such that $\lvert u^*(t)\rvert=1$ almost everywhere. Such curves are called \emph{arc length parametrized}. Due to the bracket generating condition, Rashevskii-Chow's theorem \cite[Theorem 3.31]{AgrachevBarilariBoscainBook2020} implies that the Carnot-Carth\'eodory metric given by
\begin{align*}
    d_{CC}\left( p, q \right)=\inf \{\ell(\gamma): \gamma\text{ is admissible }, \gamma(0)=p,\gamma(T)=q\
    \}
\end{align*}
is a well-defined distance function and induces the manifold topology on $M$.

An admissible curve $\gamma:[0,T]\rightarrow M$ is called a \emph{length minimizer} between $p, q \in M$ if $\gamma(0)=p$, $\gamma(T)=q$ and $\ell(\gamma)=d_{CC}\left( p, q \right)$. If $(M,d_{CC})$ is complete as a metric space, then there are length minimizers connecting any two points in $M$, see \cite[Corollary 3.49]{AgrachevBarilariBoscainBook2020}. In other words, $(M,d_{CC})$ is a \emph{geodesic complete length space} if and only if it is complete as a metric space.

Consider the cotangent bundle $T^{\ast}M$ equipped with the natural canonical form $\sigma$. The \emph{Hamiltonian vector field} of a real valued function $a \in  C^{\infty}\left(T^{\ast} M\right)$ is the unique vector field $\vec{a}$ on $T^{\ast}M$ such that
\[
    \sigma(\cdot, \vec{a}) = \mathrm{d} a,
\]
which can also be written, in the canonical coordinate chart $(x, p)$ on $T^*M$, as
\begin{align}\label{Hamiltonian vector field}
    \vec{a}=\sum_{i=1}^N\frac{\partial a}{\partial p_i}\frac{\partial}{\partial x_i}- \frac{\partial a}{\partial x_i}\frac{\partial}{\partial p_i}.
\end{align}
In the context of sub-Riemannian geometry, there is a natural smooth function on the cotangent bundle to consider. The \emph{sub-Riemannian Hamiltonian} is the map $H:T^{\ast}M\rightarrow \mathbb{R}$ defined by
\begin{equation}
    \label{eq:hamiltonian}
    H(\lambda)=\frac{1}{2}\sum_{j=1}^m h_j(\lambda)^2,
\end{equation}
where $h_j(\lambda):=\langle \lambda, Z_j(q)\rangle,$ with $q=\pi(\lambda)$. The function $H$ carries information about sub-Riemannian length minimizers through the Pontryagin Maximum Principle.
\begin{theorem}[Pontryagin Maximum Principle]\label{Pontryagin Maximum Principle}

    Suppose that $(M,\mathcal{F})$ is a sub-Riemannian structure and $\gamma:[0,T]\rightarrow M$ is a length minimizer parametrized by arc length with minimal control $u^{\ast}$. Then, there is a Lipschitz (in charts) curve $\lambda:[0, T]\rightarrow T^{\ast}M$ such that $\lambda(t)\in T^{\ast}_{\gamma(t)}M$,
    \begin{align}\label{extremal cond}
        \lambda^{\prime}(t)=\sum_{j=1}^m u^{\ast}_j(t)\vec{h}_j(\lambda(t)) \text{ a.e. } t\in [0,T]
    \end{align}
    and one of the following conditions is satisfied.
    \begin{align}\label{Pontryagin}
         & \text{ Normal: }   & h_j(\lambda(t))\equiv u^{\ast}_j(t), &  & 1 \leqslant j\leqslant m
        \\
         & \text{ Abnormal: } & h_j(\lambda(t))\equiv 0,             &  & 1 \leqslant j \leqslant m.
        \notag
    \end{align}
    Moreover, in the abnormal case, one has $\lambda(t) \neq 0$ for all $t \in [0, T]$.
\end{theorem}

The proof of this result can be found in \cite[Theorem 4.20]{AgrachevBarilariBoscainBook2020}. A curve $\lambda : [0, T] \to T^*M$ satisfying \eqref{extremal cond} is said to be a normal (resp. abnormal) extremal if it satisfies the normal (resp. abnormal) condition in \eqref{Pontryagin}. A \emph{normal} (resp. \emph{abnormal}) curve is a curve $\gamma : [0, T] \to M$ admitting a normal (resp. abnormal) lift. It is shown in \cite[Proposition 4.22]{AgrachevBarilariBoscainBook2020} that the normal condition in \eqref{Pontryagin} can be equivalently rewritten as
\begin{align}\label{e.HamiltonEquation}
    \lambda'\left( t \right)=\vec{H}\left( \lambda\left( t \right) \right).
\end{align}
Furthermore, the projection $\gamma(t) := \pi(\lambda(t))$ of a curve $\lambda : [0, T] \to T^*M$ satisfying \eqref{e.HamiltonEquation} is always a \textit{geodesic} by \cite[Theorem 4.65]{AgrachevBarilariBoscainBook2020}. We emphasize that a \emph{geodesic} is defined in this paper as a constant speed admissible curve that locally minimizes the Carnot-Carath\'{e}odory distance. Under this defintion, every constant speed length minimizer is a geodesic.

In addition to normal geodesics, a sub-Riemannian structure may admit \emph{abnormal} geodesics.
A geodesic may be both normal and abnormal, and the projection of an abnormal extremal need
not be a geodesic. Abnormal geodesics are not generally described by any ODE, and their presence can create significant
difficulties; see \cite{LudovicSR} for their role in the sub-Riemannian Monge problem.
Recent work has shown that abnormal minimizers need not be smooth
\cite{ChitourJeanMontiRiffordSacchelliSigalottiSocionovo2025}; see also
\cite{Hakavuori_2016}, where abnormal trajectories with corners are shown to be non-minimizing.
On the other hand, structures with no non-trivial abnormal trajectories are residual in the
Whitney topology \cite{ChitourJeanTrelat2006}.

A sub-Riemannian manifold that admits no non-trivial abnormal length minimizers is called \emph{ideal}.

The $\alpha$-Grushin spaces of dimension larger than 3 that we consider will \emph{not} turn out to be ideal, but this will not be an issue for us. See Theorem \ref{abnormality}.

\section{Grushin spaces and Hamiltonian theory} \label{s.GrushinSpaces}
\subsection{Grushin spaces}
Consider vector fields on $\mathbb{R}^{n+1}$ defined by
\begin{equation}
    \label{Grushin vector fields}
    \begin{split}
        X_j=      & \,\,\xi_j\partial_{x_j}                          \\
        \xi_j(x)= & \prod_{1\leqslant i\leqslant j-1}x_i^{\alpha_i},
    \end{split}
\end{equation}

for $x = (x_1, \dots, x_{n + 1})\in \mathbb{R}^{n+1}$, $1\leqslant j\leqslant n+1$ and $\alpha_i\in \mathbb{N}\cup\{0\}$ for $1\leqslant i\leqslant n$. For the rest of the paper we will use the convention that a product indexed over the empty set is identically $1$. Note that some of the vector fields vanish on the singular set $\mathcal{S}:=\{x\in \mathbb{R}^{n+1}: \prod_{j=1}^n x_j=0\}$ which is the union of the hyperplanes $\{x_j=0\}$, $1\leq j\leq n$.

The vector fields $\{X_1,..,X_{n+1}\}$ satisfy the bracket generating condition but with a varying step length. The Carnot-Carth\'{e}odory metric $d_{CC}$ is complete (by the estimate for $d_{CC}$ in \cite[Theorem 3.1]{WuJang-Mei2015}). The metric is Riemannian away from the singular set $\mathcal{S}$.
We call $\mathbb{G}^{n+1}_\alpha:=(\mathbb{R}^{n+1}, d_{CC})$ a $\alpha$-\emph{Grushin space}. Using completeness, by \cite[Theorem 3.43]{AgrachevBarilariBoscainBook2020}, there exists a length minimizer connecting any two points in the space. This turns $\mathbb{G}^{n+1}_\alpha$ into a \emph{geodesic complete length space} in the sense of \cite{BuragoBuragoIvanovBookEng}.

\subsection{Hamilton's equations for Grushin spaces}
The \emph{Hamiltonian} $H:T^{\ast}\mathbb{G}^{n+1}_\alpha\cong \mathbb{R}^{n+1}\times \mathbb{R}^{n+1}\rightarrow \mathbb{R}$   appearing in \eqref{eq:hamiltonian} is given by
\begin{equation}\label{Hamiltonian Function}
    H(x,p)=\frac{1}{2}\sum_{j=1}^{n+1}\xi_j^2(x)p_j^2 =\frac{1}{2}(p_1^2+x_1^{2\alpha_1}p_2^2+...+x_1^{2\alpha_1} ... x_n^{2\alpha_n}p_{n+1}^2).
\end{equation}
To find normal geodesics, one must solve \emph{Hamilton's equations} \eqref{e.HamiltonEquation}, which in coordinates $(x,p)$ are given by
\begin{equation}
    \label{Hamiltonian System}
    \begin{split}
        x'= & \frac{\partial H}{\partial p}
        \\
        p'= & -\frac{\partial H}{\partial x}.
    \end{split}
\end{equation}
The system \eqref{Hamiltonian System} can be expressed as a system of $2(n+1)$ ordinary differential equations.
\begin{equation}
    \label{Hamilton's Equations}
    \begin{split}
        x_j^{\prime} & = \xi_j^2 p_j,\quad 1\leqslant j\leqslant n+1
        \\
        p_j^{\prime} & = -\alpha_j \xi_j^2 x_j^{2\alpha_j-1} R_{j+1}^2,\quad 1\leqslant j\leqslant n+1
        \\
    \end{split}
\end{equation}
which are taken together with the initial conditions
\begin{align*}
    x(0)=(x_1^0,...,x_{n+1}^0)
    \\
    p(0)=(p^0_{1},...,p_{n+1}^0),
\end{align*}
and where we have recursively defined the following quantities
\begin{align*}
    R_{n+2}:= & 0,
    \\
    R_{n+1}:= & \lvert p_{n+1}\rvert,
    \\
    R_j:=     & \sqrt{R_{j+1}^2(x_j)^{2\alpha_j}+(p_j)^2},\quad   1\leqslant j\leqslant n.
\end{align*}
Our preliminary goal is to solve \eqref{Hamilton's Equations} for an initial condition $x^0\in \mathbb{G}^{n+1}_\alpha$ that is \emph{Riemannian}. A \emph{Riemannian point} $q_0\in M$ for a sub-Riemannian structure $(M,\mathcal{F})$ satisfies $\mathcal{D}_{q_0}=T_{q_0}M$. In the Grushin setting, $x^0\in \mathbb{G}^{n+1}_\alpha$ is Riemannian if $x^0_j\neq 0$ for all $j\in J$, where
\begin{align}
    J:=\{1\leq j\leq n: \alpha_j\neq 0\}.
\end{align}A point which is not Riemannian is called \emph{singular}.

For a given $x \in \mathbb{G}^{n+1}_\alpha$, we denote \begin{align}\label{H function}
    H_{x}(\cdot):=H(x,\cdot).
\end{align}
By \cite[Theorem 4.25]{AgrachevBarilariBoscainBook2020}, a geodesic $x(\cdot)=x(\,\cdot\,;p^0)$ with the initial point $x^0\in \mathbb{G}^{n+1}_\alpha$ and the initial covector $p^0$ will be arc length parametrized if and only if $p^0\in H^{-1}_{x^0}(1/2)$.  The topology of $H^{-1}_{x^0}(1/2)$ changes drastically depending on for which $i$ we have  $x^0_i=0$. We address the singular points specifically for the 3D Grushin space $\mathbb{G}^3_{(\alpha,\beta)}$ in Section ~\ref{Singular Points}.

We state the following lemma, which identifies significant constants of motion in the Hamiltonian System \eqref{Hamiltonian System}.
\begin{lemma}
    For any $j=1, ..., n+1$ we have that $R^2_j$ is constant with respect to $t$ and

    \begin{align*}
         & R^2_j=R_{j+1}^2(x^0_j)^{2\alpha_j}+(p^0_j)^2, \quad 1\leqslant j\leqslant n
        \\
         & R^2_{n+1}=(p^0_{n+1})^2.
    \end{align*}
\end{lemma}

\begin{proof}
    To see this, we start with $j=n$ and find $\frac{d}{dt}(R^2_{j+1}x^{2\alpha_j}_j+p_j^2)$ using the chain rule and \eqref{Hamilton's Equations} to see that it must be $0$. The claim follows by induction.
\end{proof}
This lemma implies that $x_1^{\prime \prime}=p_1^{\prime}=-\alpha_1x_1^{2\alpha_1-1 }R_2^2$, and $x_1^{\prime}(0)=p_1(0)=p^0_1$ and $x_1(0)=x^0_1$. This is a second order ODE whose solutions can be found explicitly using generalized trigonometric functions, which we introduce presently.

\section{Generalized trigonometric functions and solution to Hamilton's equations}\label{s.GeneralizedTrig}

\subsection{Calculus of generalized trigonometric functions} We introduce only the aspects of the generalized trigonometric functions that are necessary for our purposes, and refer to \cite{KobayashiTakeuchi2019} for a more substantial study of their properties. Note that a broader definition of generalized trigonometric functions was introduced in \cite{Lokutsievskii2019}, with applications to Hamiltonian theory. Also relevant is the paper \cite{BorzaMagnaboscoRossiTashiro2024}, which used fine properties of the functions defined in \cite{Lokutsievskii2019} to explore the MCP on sub-Finsler Heisenberg groups.

Consider for $\beta\geq 1$ a natural number, the following second order ODE.
\begin{align}\label{ODE}
    y''=        & -\beta y^{2\beta-1} \\
    \notag y(0) & =0,\,\,\,y'(0)=1.
\end{align}
We denote the solution, which is real analytic and periodic, to \eqref{ODE} by $y(\cdot)=\sin_\beta(\cdot)$. Its derivative, which is also periodic, is denoted by $y'(\cdot)=\cos_\beta(\cdot)$. The common period of $\sin_\beta,\cos_\beta$ is written as $2\pi_\beta$, and can be computed as the integral
\[\pi_\beta=2\int_0^1\frac{1}{(1-t^{2\beta})^{1/2}}\,dt.\]
A key identity throughout the paper is the energy identity relating $\sin_\beta$ and $\cos_\beta$
\begin{align}\label{gen pythag identity}
    \sin_\beta^{2\beta}+\cos_\beta^2=1.
\end{align}
Note that the general solution to the second order ODE
$y^{\prime \prime}=-\beta y^{2\beta-1}$ with $(y(0), y^{\prime}(0))\neq 0$ is given by
\begin{align}
    y(t)=A\sin_\beta(\omega t+\phi)
\end{align}
for some $A,\omega\in \mathbb{R}\setminus\{0\}$, $A\omega>0$ and $\phi\in [0,\pi_\beta)$. Whenever $\beta\neq 1$, there are no known angle addition formulae to expand $\sin_\beta(A+B)$ or $\cos_\beta(A+B)$ except for specific values of $A,B\in \mathbb{R}$.

\begin{figure}
    \centering
    \includegraphics[width=0.5\linewidth]{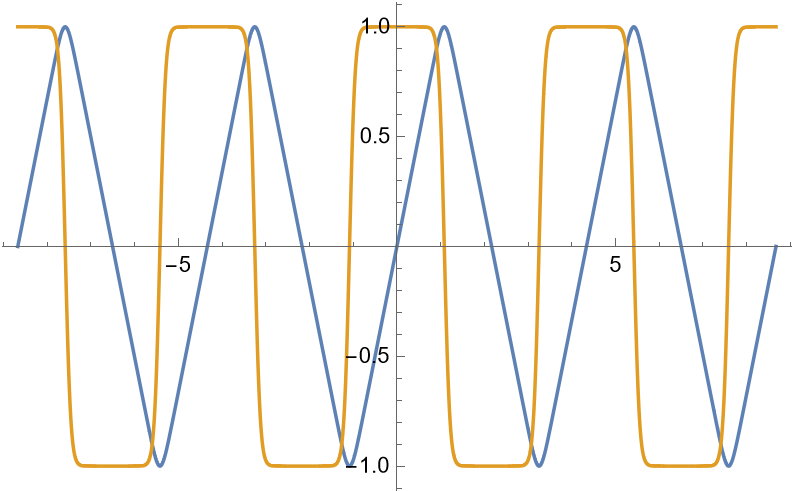}
    \caption{Graph of $y=\sin_8(x)$ in \textcolor{blue}{blue} and graph of $y=\cos_8(x)$ in \textcolor{orange}{orange}.}
    \label{fig:trig8}
\end{figure}

We record for later use that
\begin{align}\label{Taylor}
    \cos_\beta(x) & = 1-\frac{1}{2}x^{2\beta}+O(x^{2\beta+2}),\quad\quad x\to 0
    \\
    \sin_\beta(x) & = x-\frac{1}{2(2\beta+1)}x^{2\beta+1}+O(x^{2\beta+3}),\quad\quad x\to 0.
\end{align}
Note that in the case $\beta=1$, these agree with the usual power series expansions for sine and cosine. See the graph of the $\sin_\beta$ and $\cos_\beta$ functions in Figure~\ref{fig:trig8}.

Finally, we record the useful symmetry identity
\begin{align}\label{sin identity}
    \sin_{\beta}(x+\pi_\beta)=-\sin_\beta(x), x\in \mathbb{R}.
\end{align}
Then, taking derivatives, we obtain\begin{align}\label{cosine identity}
    \cos_\beta(x+\pi_\beta)=-\cos_\beta(x), x\in \mathbb{R}.
\end{align}
\subsection{Solving Hamilton's equations}\label{SHE}
We now introduce generalized spherical coordinates on \(H_{x^0}^{-1}(1/2)\), where
\(x^0\) is Riemannian. For \(\theta\geq 1\), set
\begin{equation}\label{e.rho}
    \rho_\theta(x):=|x|^{\theta-1}x .
\end{equation}
Let \(\beta=(\beta_1,\ldots,\beta_n)\), with \(\beta_j\geq 1\). Define
\(p^0=p^0(\phi)\) by
\begin{equation}\label{spherical coordinates}
    \begin{split}
        p^0_j
         & :=
        \cos_{\beta_j}(\phi_j)
        \prod_{i<j}
        \rho_{\beta_i}
        \left(
        \frac{\sin_{\beta_i}(\phi_i)}{x_i^0}
        \right),
        \qquad 1\leq j\leq n,
        \\
        p^0_{n+1}
         & :=
        \prod_{i=1}^n
        \rho_{\beta_i}
        \left(
        \frac{\sin_{\beta_i}(\phi_i)}{x_i^0}
        \right).
    \end{split}
\end{equation}
The parameter domain is
\begin{equation}\label{e.Rectangle}
    \mathcal R_\beta
    :=
    \left(\prod_{j=1}^{n-1}[0,\pi_{\beta_j}]\right)
    \times [0,2\pi_{\beta_n}).
\end{equation}
The map
\[
    \mathcal R_\beta\ni \phi
    \longmapsto
    p^0(\phi)\in H_{x^0}^{-1}(1/2)
\]
is continuous and surjective, and its restriction to
\[
    \mathcal R_\beta
    \setminus
    \left(
    \prod_{j=1}^{n-1}[0,\pi_{\beta_j}]
    \times\{0\}
    \right)
\]
is a \(C^1\)-diffeomorphism onto its image.

If \(\alpha\in(\mathbb N\cup\{0\})^n\), let
\(\widetilde\alpha\in\mathbb N^n\) be obtained by replacing each zero component of
\(\alpha\) by \(1\). For indices \(i\) with \(\alpha_i=0\), we use
\(\widetilde\alpha_i=1\) and set \(x_i^0=1\) in \eqref{spherical coordinates}. Thus
\[
    \mathcal R_\alpha:=\mathcal R_{\widetilde\alpha}
\]
parametrizes \(H_{x^0}^{-1}(1/2)\). We also set
\begin{equation}\label{e.R0}
    \mathcal R_\alpha^0
    :=
    \left\{
    \phi\in\mathcal R_{\widetilde\alpha}:
    \sin_{\widetilde\alpha_j}(\phi_j)\neq 0
    \text{ for all } j\in J
    \right\}.
\end{equation}

For \(j\in J\) and \(\phi\in\mathcal R_\alpha^0\), define
\begin{equation}\label{e.AandOmega}
    A_j
    :=
    \frac{x_j^0}{\sin_{\alpha_j}(\phi_j)},
    \qquad
    \omega_j
    :=
    \frac{\sin_{\alpha_j}(\phi_j)}{x_j^0}
    \prod_{i<j}
    \rho_{\alpha_i}
    \left(
    \frac{\sin_{\alpha_i}(\phi_i)}{x_i^0}
    \right).
\end{equation}
It follows that
\begin{align}
    A_j\omega_j\cos_{\alpha_j}(\phi_j)
     & =
    p_j^0,
    \label{properties}
    \\
    \frac{\omega_j^2}{A_j^{2(\alpha_j-1)}}
     & =
    R_{j+1}^2.
    \label{properties 2}
\end{align}
Moreover, \(R_{j+1}\neq 0\) if and only if
\(\sin_{\alpha_i}(\phi_i)\neq 0\) for every \(i\in J\) with \(i\leq j\).
\begin{definition}
    Let \(x(t)\) be a normal geodesic in \(\mathbb G_\alpha^{n+1}\) with initial point
    \(x^0\). We say that \(x(t)\) is \emph{non-trivial} if it is not strictly contained in any
    hyperplane \(x_j=x_j^0\), \(j\in J\).
\end{definition}
Note that
non-trivial geodesics correspond precisely to \(\phi\in\mathcal R_\alpha^0\). Define
\begin{equation}\label{x definition}
    f_j(t;x^0,p^0,\alpha)
    =
    \begin{cases}
        x_j^0+p_j^0t,
         & R_{j+1}=0 \text{ or } \alpha_j=0,          \\[2mm]
        A_j\sin_{\alpha_j}(\omega_jt+\phi_j),
         & R_{j+1}\neq 0 \text{ and } \alpha_j\geq 1,
    \end{cases}
\end{equation}
and
\begin{equation}\label{p definition}
    g_j(t;x^0,p^0,\alpha)
    =
    \begin{cases}
        p_j^0,
         & R_{j+1}=0 \text{ or } \alpha_j=0,          \\[2mm]
        A_j\omega_j\cos_{\alpha_j}(\omega_jt+\phi_j),
         & R_{j+1}\neq 0 \text{ and } \alpha_j\geq 1.
    \end{cases}
\end{equation}
Set \(x_1(t)=f_1(t;x^0,p^0,\alpha)\) and \(p_1(t)=g_1(t;x^0,p^0,\alpha)\). Inductively, for
\(2\leq j\leq n+1\), define
\begin{equation}\label{e.inductive-coordinates}
    x_j(t)
    :=
    f_j\left(\int_0^t \xi_j^2(s)\,ds;x^0,p^0,\alpha\right),
    \qquad
    p_j(t)
    :=
    g_j\left(\int_0^t \xi_j^2(s)\,ds;x^0,p^0,\alpha\right).
\end{equation}
The initial conditions are immediate. In the non-linear case, using
\[
    \frac{\omega_j^2}{A_j^{2(\alpha_j-1)}}=R_{j+1}^2
\]
gives
\begin{align*}
    \frac{d}{dt}
    \left[
        A_j\omega_j
        \cos_{\alpha_j}
        \left(
        \omega_j\int_0^t \xi_j^2(s)\,ds+\phi_j
        \right)
        \right]
     & =
    -\alpha_j\xi_j^2(t)A_j\omega_j^2
    \sin_{\alpha_j}^{2\alpha_j-1}
    \left(
    \omega_j\int_0^t \xi_j^2(s)\,ds+\phi_j
    \right)
    \\
     & =
    -\alpha_j\xi_j^2(t)x_j(t)^{2\alpha_j-1}R_{j+1}^2.
\end{align*}
Thus \((x_j,p_j)\) satisfies Hamilton's equations.

\begin{theorem}\label{Solution to Hamilton}
    Let \(\alpha\in(\mathbb N\cup\{0\})^n\), and let
    \(\mathbb G_\alpha^{n+1}=(\mathbb R^{n+1},d_{CC})\). Every arclength-parametrized normal
    geodesic \(x:[0,T]\to\mathbb G_\alpha^{n+1}\) with Riemannian initial point \(x^0\) has a
    unique lift \(\lambda=(x,p)\) solving Hamilton's equations with
    \(p(0)=p^0\in H_{x^0}^{-1}(1/2)\). Its coordinates are
    \begin{equation}\label{e.coordinates}
        x_j(t)
        =
        f_j\left(\int_0^t \xi_j^2(s)\,ds;x^0,p^0,\alpha\right),
        \qquad
        p_j(t)
        =
        g_j\left(\int_0^t \xi_j^2(s)\,ds;x^0,p^0,\alpha\right),
    \end{equation}
    where
    \[
        \xi_j(t)=\prod_{i<j}x_i(t)^{\alpha_i}.
    \]
\end{theorem}

\begin{remark}
    With the convention \(\zeta^{2a}:=(\zeta^2)^a\), the same formulas hold for
    \(\alpha_j\in[1,\infty)\cup\{0\}\) if
    \[
        X_1=\partial_{x_1},
        \qquad
        X_j=|x_1|^{\alpha_1}\cdots |x_{j-1}|^{\alpha_{j-1}}\partial_{x_j},
        \quad 2\leq j\leq n+1.
    \]
    These length spaces, studied in \cite{WuJang-Mei2015, Borza2022}, need not be smooth
    sub-Riemannian manifolds in the sense of \cite{AgrachevBarilariBoscainBook2020}, so we restrict to integer \(\alpha_j\).
\end{remark}

We next record a recursion formula. For \(\beta\geq 1\), set
\[
    \eta_\beta(s):=s-\sin_\beta(s)\cos_\beta(s).
\]

\begin{proposition}\label{Recursion Result}
    If \(j\in J\) and \(R_{j+1}>0\), then
    \begin{equation}\label{e.recursion-integral}
        \int_0^t \xi_{j+1}^2(s)\,ds
        =
        \frac{A_j^{2\alpha_j}}{\omega_j(\alpha_j+1)}
        \left[
            \eta_{\alpha_j}
            \left(
            \omega_j\int_0^t \xi_j^2(s)\,ds+\phi_j
            \right)
            -
            \eta_{\alpha_j}(\phi_j)
            \right].
    \end{equation}
\end{proposition}

\begin{proof}
    Since
    \[
        \eta_\beta'(s)=(\beta+1)\sin_\beta^{2\beta}(s),
    \]
    the result follows from the change of variables
    \[
        r=\omega_j\int_0^s \xi_j^2(u)\,du+\phi_j
    \]
    in the integral
    \[
        \int_0^t \xi_{j+1}^2(s)\,ds
        =
        A_j^{2\alpha_j}
        \int_0^t
        \xi_j^2(s)
        \sin_{\alpha_j}^{2\alpha_j}
        \left(
        \omega_j\int_0^s \xi_j^2(u)\,du+\phi_j
        \right)
        ds.
    \]
\end{proof}

Consequently, if \(R_{j+2}>0\) and \(\alpha_{j+1}\geq 1\), then
\begin{equation}\label{Recursion 1}
    x_{j+1}(t)
    =
    A_{j+1}
    \sin_{\alpha_{j+1}}
    \left(
    \frac{\omega_{j+1}A_j^{2\alpha_j}}{\omega_j(\alpha_j+1)}
    \left[
        \eta_{\alpha_j}
        \left(
        \omega_j\int_0^t \xi_j^2(s)\,ds+\phi_j
        \right)
        -
        \eta_{\alpha_j}(\phi_j)
        \right]
    +\phi_{j+1}
    \right).
\end{equation}
If \(R_{j+2}=0\) or \(\alpha_{j+1}=0\), then
\begin{equation}\label{Recursion 2}
    x_{j+1}(t)
    =
    x_{j+1}^0
    +
    \frac{p_{j+1}^0A_j^{2\alpha_j}}{\omega_j(\alpha_j+1)}
    \left[
        \eta_{\alpha_j}
        \left(
        \omega_j\int_0^t \xi_j^2(s)\,ds+\phi_j
        \right)
        -
        \eta_{\alpha_j}(\phi_j)
        \right].
\end{equation}
In particular, \eqref{Recursion 2} applies when \(j=n\). Equations
\eqref{Recursion 1}--\eqref{Recursion 2} express the geodesics recursively in generalized
trigonometric functions.
\section{Abnormal minimizers in Grushin spaces}

Although \(\mathbb G^2_\alpha\) is ideal for \(\alpha\geq 1\), higher-dimensional
Grushin spaces need not be. The obstruction comes from lower-dimensional Grushin
geodesics contained in singular hyperplanes. We show that these are the only abnormal
minimizers, and that none of them are strictly abnormal.

\begin{theorem}\label{3D abnormal}
    Let \(\alpha,\beta\in \mathbb{N}\) In \(\mathbb G^3_{(\alpha,\beta)}\), the only non-trivial
    abnormal length minimizers starting at \(q_0=(x_0,y_0,z_0)\) occur when \(y_0=0\).
    They are, up to reparametrization,
    \[
        \gamma(t)=(x_0+\eta t,0,z_0),
        \qquad \eta\neq 0,
    \]
    and are also normal.
\end{theorem}

\begin{proof}
    Write the generating vector fields of $\mathbb{G}_{(\alpha,\beta)}$ as
    \[
        X=\partial_x,\qquad
        Y=x^\alpha\partial_y,\qquad
        Z=x^\alpha y^\beta\partial_z,
    \]
    and let \(\lambda=(\gamma,p)\), \(p=(u,v,w)\), be an abnormal lift of a minimizing
    trajectory \(\gamma=(x,y,z)\). The abnormal conditions are
    \[
        h_X=u=0,\qquad
        h_Y=x^\alpha v=0,\qquad
        h_Z=x^\alpha y^\beta w=0,
    \]
    with \(p\not\equiv 0\). The adjoint equations give \(\dot w=0\). If \(v\) were nonzero
    on some interval, then \(x\equiv 0\) there, forcing \(\gamma\) to be constant on a
    subinterval, contradicting minimality. Hence \(v\equiv 0\), so \(w\equiv w_0\neq 0\).
    Then \(h_Z=0\) forces \(y\equiv 0\); otherwise the same argument would again force a
    constant subarc. Thus \(z\equiv z_0\), and \(\gamma\) lies in the line
    \(\{(x,0,z_0):x\in\mathbb R\}\).

    For its minimal control \(\mu=(\eta^*,\epsilon^*,\delta^*)\), admissibility along this
    line gives \(\epsilon^*=0\), and
    \[
        \ell(\gamma)
        =
        \int_0^T\sqrt{(\eta^*)^2+(\delta^*)^2}\,dt
        \geq
        \left|\int_0^T\eta^*(t)\,dt\right|
        =
        |x(T)-x(0)|.
    \]
    Equality is necessary for minimality, so \(\delta^*=0\) and \(\eta^*=\pm1\) after
    arclength reparametrization. The resulting line has the normal covector lift
    \(p\equiv(\pm1,0,0)\).
\end{proof}

The following elementary projection argument records the lower-triangular structure of the
Grushin vector fields.

\begin{theorem}\label{isometric embeddings}
    Let \(n\geq 1\), \(1\leq k\leq n\), and
    \[
        \alpha^k:=(\alpha_1,\ldots,\alpha_{n-k}),
    \]
    with \(\mathbb G^1_{()}:=\mathbb R\). Suppose \(q_0,q_1\in\mathbb G^{n+1}_\alpha\)
    agree in their last \(k\) coordinates, and write \(q_0^k,q_1^k\) for their projections
    onto the first \(n+1-k\) coordinates. Then
    \[
        d_{\mathbb G^{n+1}_\alpha}(q_0,q_1)
        =
        d_{\mathbb G^{n+1-k}_{\alpha^k}}(q_0^k,q_1^k).
    \]
    Moreover, a curve \(\gamma\) from \(q_0\) to \(q_1\) is minimizing in
    \(\mathbb G^{n+1}_\alpha\) with the last \(k\) coordinates fixed if and only if its
    projection \(\gamma^k\) is minimizing in \(\mathbb G^{n+1-k}_{\alpha^k}\). Hence
    \[
        \mathbb R
        \hookrightarrow
        \mathbb G^2_{\alpha^{n-1}}
        \hookrightarrow
        \cdots
        \hookrightarrow
        \mathbb G^{n+1}_\alpha
    \]
    is a chain of isometric embeddings.
\end{theorem}

\begin{proof}
    Let \(\gamma=(x_1,\ldots,x_{n+1})\) be admissible with minimal control
    \(u^*=(u_1^*,\ldots,u_{n+1}^*)\). Its projection
    \[
        \gamma^k=(x_1,\ldots,x_{n+1-k})
    \]
    is admissible in \(\mathbb G^{n+1-k}_{\alpha^k}\), with control
    \((u_1^*,\ldots,u_{n+1-k}^*)\). Therefore
    \[
        \ell_{\mathbb G^{n+1}_\alpha}(\gamma)
        =
        \int_0^T
        \sqrt{
            \sum_{j=1}^{n+1-k}(u_j^*)^2
            +
            \sum_{j=n+2-k}^{n+1}(u_j^*)^2
        }\,dt
        \geq
        \ell_{\mathbb G^{n+1-k}_{\alpha^k}}(\gamma^k).
    \]
    Taking the infimum over competitors gives
    \[
        d_{\mathbb G^{n+1}_\alpha}(q_0,q_1)
        \geq
        d_{\mathbb G^{n+1-k}_{\alpha^k}}(q_0^k,q_1^k).
    \]
    The reverse inequality follows by lifting any admissible curve in
    \(\mathbb G^{n+1-k}_{\alpha^k}\) while keeping the last \(k\) coordinates fixed. The
    statements about minimizers and the isometric embeddings follow immediately.
\end{proof}

\begin{theorem}\label{abnormality}
    Let \(n\geq 2\) and \(\alpha_j\geq 1\) for all \(j\). A non-trivial abnormal length
    minimizer in \(\mathbb G^{n+1}_\alpha\) can occur only when the initial point lies on one
    of the hyperplanes
    \[
        x_j=0,\qquad 2\leq j\leq n.
    \]
    More precisely, if \(j=n+2-k\) with \(2\leq k\leq n\), then the abnormal minimizers
    contained in \(x_j=0\) are exactly the curves
    \[
        \gamma(t)
        =
        \bigl(\gamma^k(t),0,x_{j+1}^0,\ldots,x_{n+1}^0\bigr),
    \]
    where \(\gamma^k\) is a length minimizer in
    \(\mathbb G^{n+1-k}_{\alpha^k}\). In particular, every abnormal length minimizer is also
    normal.
\end{theorem}

\begin{proof}
    Let \(\lambda=(x,p)\) be an abnormal lift of a unit-speed minimizing curve
    \[x(t)=(x_1(t),\ldots,x_{n+1}(t)).\] Since
    \[
        h_i(\lambda)=\xi_i(x)p_i=0,
        \qquad 1\leq i\leq n+1,
    \]
    and \(p\not\equiv0\), let \(K\geq 3\) be the first index with \(p_K\not\equiv0\). Then
    \[
        0=h_K(\lambda)=\xi_K(x)p_K
    \]
    forces \(x_j=0\) on some interval for some \(2\leq j<K\). By admissibility, all
    coordinates \(x_{j+1},\ldots,x_{n+1}\) are constant on that interval. Maximality of the
    interval and continuity give \(x_j\equiv0\) on all of \([0,T]\); otherwise, immediately
    after an endpoint of the maximal interval, some earlier coordinate must vanish and again
    force \(x_j\) to remain constant, contradicting maximality. Hence
    \(x_{j+1},\ldots,x_{n+1}\) are constant on \([0,T]\).

    Writing \(j=n+2-k\), Theorem~\ref{isometric embeddings} shows that the remaining
    coordinates form a minimizer \(\gamma^k\) in
    \(\mathbb G^{n+1-k}_{\alpha^k}\). Conversely, any such lifted minimizer is abnormal: a
    constant covector supported in the frozen coordinates annihilates the distribution because
    \(x_j\equiv0\). It is also normal by lifting a normal extremal for \(\gamma^k\) and adding
    zero covector components in the frozen directions. Thus no abnormal minimizer is strictly
    abnormal.
\end{proof}
Consequently, the only failure of ideality comes from lower-dimensional minimizing normal geodesics embedded in singular hyperplanes. Away from the hyperplanes \(x_j=0\),
\(2\leq j\leq n\), all length minimizers are normal geodesics.
\section{Cut and conjugate times in Grushin spaces}
\label{s.CutConjGrushin}
For a geodesic (either normal or abnormal) \(\gamma(t)\) starting from \(q_0\), define
\[
    t_{\operatorname{cut}}(\gamma)
    :=
    \sup\{t>0:\gamma|_{[0,t]}\text{ is minimizing}\}.
\]
When $\gamma(t)=\gamma(t;p^0)$ is normal, the cut time may be written as a function of the initial covector. For a normal geodesic, which may also be abnormal, we define
\[
    t_{\operatorname{con}}(\gamma)
    :=
    \inf\{t>0:\exp_{q_0}\text{ is singular at }tp^0\}.
\]
Here \(\exp_{q_0}(p^0)=\gamma(1;p^0)\), and by homogeneity
\[
    \exp_{q_0}(tp^0)=\gamma(t;p^0).
\]
We compute conjugate times by writing \(p^0=p^0(\phi)\) in generalized spherical
coordinates and evaluating the Jacobian determinant
\[
    D(t,\phi)
    :=
    \det \frac{\partial \exp_{x^0}(tp^0(\phi))}{\partial(t,\phi_1,\ldots,\phi_n)}.
\]
Let \(x^0\in\mathbb G^{n+1}_\alpha\) be Riemannian and
\(\phi\in\mathcal R_\alpha^0\). Set
\begin{equation}\label{delta}
    \delta_j
    :=
    \begin{cases}
        \rho_{\alpha_j}\left(\dfrac{\sin_{\alpha_j}(\phi_j)}{x_j^0}\right),
         & j\in J,    \\[2mm]
        \sin(\phi_j),
         & j\notin J,
    \end{cases}
    \qquad 1\leq j\leq n.
\end{equation}
Then
\[
    p_j^0
    =
    \begin{cases}
        \cos_{\alpha_j}(\phi_j)\prod_{i<j}\delta_i,
         & j\in J,              \\[2mm]
        \cos(\phi_j)\prod_{i<j}\delta_i,
         & j\notin J,\ j\leq n, \\[2mm]
        \prod_{i=1}^n\delta_i,
         & j=n+1.
    \end{cases}
\]
For \(j\in J\),
\[
    \omega_j
    =
    \frac{\sin_{\alpha_j}(\phi_j)}{x_j^0}
    \prod_{i<j}\delta_i,
    \qquad
    A_j
    =
    \frac{x_j^0}{\sin_{\alpha_j}(\phi_j)}.
\]
Define
\begin{equation}\label{e.Qj}
    Q_j(t,\phi)
    :=
    \omega_j\int_0^t\xi_j^2(s)\,ds+\phi_j,
\end{equation}
\[
    S_j(t,\phi)
    :=
    \begin{cases}
        \displaystyle \int_0^t\xi_j^2(s)\,ds,
         & j\notin J, \\[3mm]
        \displaystyle
        \cos_{\alpha_j}(Q_j)
        \left(
        \cot_{\alpha_j}(\phi_j)\omega_j
        \int_0^t\xi_j^2(s)\,ds
        +1
        \right)
        -
        \cot_{\alpha_j}(\phi_j)\sin_{\alpha_j}(Q_j),
         & j\in J.
    \end{cases}
\]
Finally, set
\begin{equation}\label{e.gphi}
    g(\phi)
    :=
        (-1)^{|J|}
    \left(
    \prod_{\substack{1\leq j\leq n+1\\ j-1\notin J}}
    \prod_{i<j}\delta_i
    \right)
    \prod_{j\in J}
    \left|
    \frac{x_j^0}{\sin_{\alpha_j}(\phi_j)}
    \right|^{\alpha_j+1}.
\end{equation}

\begin{lemma}\label{Determinant Formula Theorem}
    For \(n\geq1\), \(\alpha\in(\mathbb N\cup\{0\})^n\), and
    \(\phi\in\mathcal R_\alpha^0\), the Jacobian determinant of the exponential map in
    coordinates \((t,\phi)\) satisfies
    \begin{equation}\label{Determinant Factorization}
        D(t,\phi)=g(\phi)\prod_{j=1}^n S_j(t,\phi).
    \end{equation}
\end{lemma}

\begin{proof}
    We give the proof when \(\alpha_j\geq1\) for all \(j\); the case with zero components is
    identical after replacing the corresponding factors by the Euclidean ones.

    The case \(n=1\) is the determinant computation for the \(\alpha\)-Grushin plane and may be found in
    \cite{Borza2022}. Assume the formula holds in dimension \(n+1\), and put
    \[
        \widetilde\alpha=(\alpha_1,\ldots,\alpha_{n+1}).
    \]
    Let
    \[
        \widetilde x(t,\widetilde\phi)
        =
        (\widetilde x_1,\ldots,\widetilde x_{n+2})
    \]
    be the geodesic in \(\mathbb G^{n+2}_{\widetilde\alpha}\), with
    \(\widetilde\phi=(\phi_1,\ldots,\phi_{n+1})\). Since
    \(\widetilde x_j=x_j\) for \(1\leq j\leq n\), expanding the determinant along the last
    parameter gives
    \begin{align}\label{Determinant Induction}
        \widetilde D
         & =
        \partial_{\phi_{n+1}}\widetilde x_{n+2}
        \sum_{\sigma\in S_{n+1}}
        \operatorname{sgn}(\sigma)
        \partial_{\sigma(t)}x_1\cdots
        \partial_{\sigma(\phi_n)}\widetilde x_{n+1}
        \notag   \\
         & \quad
        -
        \partial_{\phi_{n+1}}\widetilde x_{n+1}
        \sum_{\sigma\in S_{n+1}}
        \operatorname{sgn}(\sigma)
        \partial_{\sigma(t)}x_1\cdots
        \partial_{\sigma(\phi_n)}\widetilde x_{n+2}.
    \end{align}
    Using the recursion formula for the last two coordinates and the identity \eqref{gen pythag identity} one obtains
    \[
        \partial_{\phi_{n+1}}\widetilde x_{n+2}
        =
        -
        \cos_{\alpha_{n+1}}(Q_{n+1})
        S_{n+1}
        \left|
        \frac{x_{n+1}^0}{\sin_{\alpha_{n+1}}(\phi_{n+1})}
        \right|^{\alpha_{n+1}+1},
    \]
    and
    \[
        \partial_{\phi_{n+1}}\widetilde x_{n+1}
        =
        \frac{x_{n+1}^0}{\sin_{\alpha_{n+1}}(\phi_{n+1})}
        S_{n+1}.
    \]
    The same recursion gives, for derivatives with respect to \(t,\phi_1,\ldots,\phi_n\),
    \[
        \partial \widetilde x_{n+1}
        =
        \cos_{\alpha_{n+1}}(Q_{n+1})\,\partial x_{n+1},
    \]
    and
    \[
        \partial \widetilde x_{n+2}
        =
        \frac{
            \sin_{\alpha_{n+1}}(\phi_{n+1})
            |x_{n+1}^0|^{\alpha_{n+1}+1}
        }{
            x_{n+1}^0
            |\sin_{\alpha_{n+1}}(\phi_{n+1})|^{\alpha_{n+1}+1}
        }
        \sin_{\alpha_{n+1}}^{2\alpha_{n+1}}(Q_{n+1})
        \,\partial x_{n+1}.
    \]
    Substitution into \eqref{Determinant Induction} yields
    \[
        \widetilde D
        =
        -
        \left|
        \frac{x_{n+1}^0}{\sin_{\alpha_{n+1}}(\phi_{n+1})}
        \right|^{\alpha_{n+1}+1}
        S_{n+1}
        \left(
        \cos_{\alpha_{n+1}}^2(Q_{n+1})
        +
        \sin_{\alpha_{n+1}}^{2\alpha_{n+1}}(Q_{n+1})
        \right)
        D.
    \]
    The generalized Pythagorean identity reduces this to
    \[
        \widetilde D
        =
        -
        \left|
        \frac{x_{n+1}^0}{\sin_{\alpha_{n+1}}(\phi_{n+1})}
        \right|^{\alpha_{n+1}+1}
        S_{n+1}D.
    \]
    Applying the induction hypothesis gives
    \[
        \widetilde D
        =
        \widetilde g(\widetilde\phi)
        \prod_{j=1}^{n+1}S_j,
    \]
    as required.
\end{proof}

Since \(g(\phi)\neq0\) on \(\mathcal R_\alpha^0\), the first conjugate time is determined
by the first zero of one of the factors \(S_j\):
\[
    t_{\operatorname{con}}(\phi)
    =
    \min_{1\leq j\leq n}
    \inf\{t>0:S_j(t,\phi)=0\}.
\]
If \(j\notin J\), then \(S_j(t,\phi)=\int_0^t\xi_j^2(s)\,ds>0\) for \(t>0\). Hence only
indices \(j\in J\) matter for conjugacy. For such \(j\), define
\begin{equation}\label{Conjectured Cut}
    \tau_j(\phi)
    :=
    \min
    \left\{
    t>0:
    |\omega_j|\int_0^t\xi_j^2(s)\,ds
    =
    \pi_{\alpha_j}
    \right\}.
\end{equation}

\begin{proposition}\label{tan function solution}
    Let \(\alpha\geq1\) and
    \[
        \phi\in[0,2\pi_\alpha)
        \setminus
        \left\{
        \frac{\pi_\alpha}{2},
        \frac{3\pi_\alpha}{2}
        \right\}.
    \]
    Then the equation
    \begin{equation}\label{tan equation}
        \tan_\alpha(\phi+\zeta)-\tan_\alpha(\phi)=\zeta
    \end{equation}
    has no nonzero solution with \(|\zeta|<\pi_\alpha\).
\end{proposition}

\begin{proof}
    On any interval on which \(\tan_\alpha(\phi+\zeta)\) is finite,
    \[
        \frac{d}{d\zeta}\tan_\alpha(\phi+\zeta)
        =
        \frac{
            1+(\alpha-1)\sin_\alpha^{2\alpha}(\phi+\zeta)
        }{
            \cos_\alpha^2(\phi+\zeta)
        }
        \geq 1.
    \]
    Thus
    \[
        \zeta\mapsto \tan_\alpha(\phi+\zeta)-\tan_\alpha(\phi)-\zeta
    \]
    is monotone on each such interval and vanishes at \(\zeta=0\). The vertical asymptotes of
    \(\tan_\alpha\) separate the remaining subintervals of \((-\pi_\alpha,\pi_\alpha)\), and on
    each of them the sign is fixed by the one-sided limits. Hence no other zero occurs.
\end{proof}

\begin{theorem}\label{Determinant}
    Let \(x^0\in\mathbb G_\alpha^{n+1}\) be Riemannian, and let
    \(\phi\in\mathcal R_\alpha^0\) correspond to a non-trivial geodesic. Define
    \[
        \tau(\phi):=\min_{j\in J}\tau_j(\phi),
    \]
    with \(\tau_j\) as in \eqref{Conjectured Cut}. Then
    \begin{equation}\label{conjugate inequality}
        0<\tau(\phi)\leq t_{\operatorname{con}}(\phi).
    \end{equation}
    Moreover, equality holds if and only if there is \(j\in J\) such that
    \(\tau=\tau_j\) and
    \[
        \phi_j=\frac{\pi_{\alpha_j}}{2}
        \quad\text{or}\quad
        \phi_j=\frac{3\pi_{\alpha_j}}{2}.
    \]
\end{theorem}

\begin{proof}
    For \(j\notin J\), \(S_j\) has no positive zero. Let \(j\in J\). If
    \[
        \phi_j=\frac{\pi_{\alpha_j}}{2}
        \quad\text{or}\quad
        \phi_j=\frac{3\pi_{\alpha_j}}{2},
    \]
    then \(S_j(t,\phi)=0\) first occurs precisely when
    \[
        |\omega_j|\int_0^t\xi_j^2(s)\,ds=\pi_{\alpha_j},
    \]
    that is, at \(t=\tau_j\).

    Otherwise, for \(0<t<\tau_j\), the equation \(S_j(t,\phi)=0\) is equivalent to
    \[
        \tan_{\alpha_j}
        \left(
        \phi_j+\omega_j\int_0^t\xi_j^2(s)\,ds
        \right)
        -
        \tan_{\alpha_j}(\phi_j)
        =
        \omega_j\int_0^t\xi_j^2(s)\,ds.
    \]
    By Proposition~\ref{tan function solution}, this has no solution for \(0<t<\tau_j\).
    Thus no factor \(S_j\) vanishes before \(\tau=\min_{j\in J}\tau_j\), proving
    \(\tau\leq t_{\operatorname{con}}\). The equality statement follows from the special
    cases above.
\end{proof}

We next obtain the corresponding upper bound for the cut time.

\begin{theorem}\label{upper bound cut time}
    Let \(x(t;\phi)\) be a non-trivial geodesic in \(\mathbb G^{n+1}_\alpha\). Then
    \begin{equation}\label{cut time inequality}
        t_{\operatorname{cut}}(\phi)\leq \tau(\phi).
    \end{equation}
\end{theorem}

\begin{proof}
    Let \(\tau(\phi)=\tau_j(\phi)\) for some \(j\in J\), and define \(\phi'\) by
    \[
        \phi_i'=\phi_i\quad(i\neq j),
        \qquad
        \phi_j'=\pi_{\alpha_j}-\phi_j.
    \]
    Then \(|\omega_j(\phi')|=|\omega_j(\phi)|\), so
    \[
        \tau_j(\phi')=\tau_j(\phi).
    \]
    Moreover, for \(i<j\), the coordinate \(x_i(t;\phi)\) is independent of \(\phi_j\), hence
    \[
        x_i(\tau_j;\phi')=x_i(\tau_j;\phi),
        \qquad i<j.
    \]
    At \(t=\tau_j\), the symmetry identities \eqref{sin identity} and \eqref{cosine identity} for \(\sin_{\alpha_j}\) and
    \(\cos_{\alpha_j}\) give
    \[
        x_j(\tau_j;\phi')=x_j(\tau_j;\phi)=-x_j^0.
    \]
    Finally, the recursive formulae \eqref{Recursion 1}--\eqref{Recursion 2} imply the same
    equality for all later coordinates:
    \[
        x_i(\tau_j;\phi')=x_i(\tau_j;\phi),
        \qquad j<i\leq n+1.
    \]
    Thus the two distinct geodesics \(x(t;\phi)\) and \(x(t;\phi')\) meet at time \(\tau_j\).
    Therefore \(x(t;\phi)\) cannot minimize beyond \(\tau_j=\tau(\phi)\).
\end{proof}

The preceding argument gives only an upper bound: it does not exclude the possibility that
\(x(t;\phi)\) meets a third minimizing geodesic before \(\tau(\phi)\). We conjecture that
this never occurs.

\begin{conj}\label{conjecture}
    For every non-trivial geodesic in \(\mathbb G^{n+1}_\alpha\),
    \[
        t_{\operatorname{cut}}(\phi)
        =
        \tau(\phi)
        =
        \min_{j\in J}\tau_j(\phi),
    \]
    where \(\tau_j\) is defined in \eqref{Conjectured Cut}.
\end{conj}

For \(n=1\), this is the known cut-time formula for the \(\alpha\)-Grushin plane
\cite{AgrachevBarilariBoscainBook2020,Borza2022}. We prove the conjecture for \(n=2\) in
the next two sections.

The standard approach for studying cut times in sub-Riemannian manifolds is the extended Hadamard technique, which we recall here and whose proof can be found in \cite[Chapter 13]{AgrachevBarilariBoscainBook2020}.
\begin{theorem}[Section~13.4 in \cite{AgrachevBarilariBoscainBook2020}]\label{Extended Hadamard Technique} Let $M$ be a complete sub-Riemannian manifold and let $q_0\in M$ be a Riemannian point. Denote by $\operatorname{Cut}^{\ast}(q_0) \subset M$ the conjectured cut locus of $q_0$, and let $\tau^{\ast}(\lambda_0) \in (0, +\infty]$ be the conjectured cut time for an arc-length geodesic with initial covector $\lambda_0 \in H^{-1}_{q_0}(1/2)$.

    Define $N$ to be the conjectured cotangent injectivity domain, that is, the set of covectors in $T^{\ast}_{q_0}M$ for which the associated geodesics are conjectured to remain optimal up to time $1$:
    \[
        N := \{t\lambda_0 : \lambda_0 \in H^{-1}_{q_0}(1/2),\ t \in [0, \tau^{\ast}(\lambda_0))\}.
    \]
    Assume that $N$ satisfies the following properties.
    \begin{enumerate}[label=\normalfont\roman*), itemsep=4pt]
        \item $\exp_{q_0}(N) = M \setminus \operatorname{Cut}^{\ast}(q_0)$;
        \item The restriction $\left.\exp_{q_0}\right|_N$ is a proper map, and is invertible at every point of $N$;
        \item The set $M \setminus \operatorname{Cut}^\ast(q_0)$ is simply connected.
    \end{enumerate}
    Then the conjectured cut time and cut locus are exact. In other words, $t_{\operatorname{cut}} = \tau^{\ast}$ and $\operatorname{Cut}(q_0) = \operatorname{Cut}^{\ast}(q_0)$.
\end{theorem}
We carry out the necessary steps to apply this theorem to the Grushin spaces with $n=2$ in Theorem \ref{main theorem}. The main challenge is to verify assumptions i) and iii) beyond $n=2$. However, Theorem~\ref{Determinant} allows us to prove ii) for any $n$.

\subsection{Conjugate times in \texorpdfstring{$\mathbb G^3_{(\alpha,\beta)}$}{G3}}

Let \(\phi=(\phi_1,\phi_2)\) be the generalized spherical coordinates on
\(H_{q_0}^{-1}(1/2)\). For non-trivial geodesics,
\[
    \phi\in \mathcal R^0_{(\alpha,\beta)}
    =
    (0,\pi_\alpha)\times
    \bigl((0,2\pi_\beta)\setminus\{0,\pi_\beta\}\bigr).
\]
Write
\[
    A:=Q_1=\omega_1t+\phi_1,
    \qquad
    B:=Q_2=\omega_2\int_0^t x^{2\alpha}(s)\,ds+\phi_2.
\]
The determinant
\[
    D(t,\phi_1,\phi_2)
    =
    \det
    \frac{\partial(x,y,z)}{\partial(t,\phi_1,\phi_2)}
\]
is, by Lemma~\ref{Determinant Formula Theorem},
\begin{equation}\label{det formula 3D}
    \begin{aligned}
        D(t,\phi_1,\phi_2)
         & =
        -
        \frac{|x_0|^{\alpha+1}|y_0|^{\beta+1}}
        {|\sin_\alpha(\phi_1)|^{\alpha+1}|\sin_\beta(\phi_2)|^{\beta+1}}
        \\
         & \quad\times
        \left[
            \cos_\alpha(A)
            \left(\cot_\alpha(\phi_1)\omega_1t+1\right)
            -
            \cot_\alpha(\phi_1)\sin_\alpha(A)
            \right]
        \\
         & \quad\times
        \left[
            \cos_\beta(B)
            \left(
            \cot_\beta(\phi_2)\omega_2
            \int_0^t x^{2\alpha}(s)\,ds
            +1
            \right)
            -
            \cot_\beta(\phi_2)\sin_\beta(B)
            \right].
    \end{aligned}
\end{equation}
For \(\lambda_0=(u_0,v_0,w_0)\in H_{q_0}^{-1}(1/2)\), define
\[
    \tau_1(\lambda_0)
    :=
    \begin{cases}
        \dfrac{\pi_\alpha}{|\omega_1|},
         & (u_0,v_0)\neq(\pm1,0), \\[2mm]
        +\infty,
         & (u_0,v_0)=(\pm1,0),
    \end{cases}
\]
\[
    \tau_2(\lambda_0)
    :=
    \begin{cases}
        \min\left\{
        t>0:
        |\omega_2|\displaystyle\int_0^t x^{2\alpha}(s)\,ds=\pi_\beta
        \right\},
         & w_0\neq0, \\[3mm]
        +\infty,
         & w_0=0.
    \end{cases}
\]
Set \(\tau=\min\{\tau_1,\tau_2\}\).

\begin{lemma}\label{Con Times 3D}
    Let \(\gamma(t)=\exp_{q_0}(t\lambda_0)\) be a geodesic in
    \(\mathbb G^3_{(\alpha,\beta)}\) starting from a Riemannian point. Then \(\gamma\) has no
    conjugate time in \((0,\tau)\). Moreover, \(t_{\operatorname{con}}=\tau\) if either
    \[
        \tau=\tau_1
        \quad\text{and}\quad
        \phi_1=\frac{\pi_\alpha}{2},
    \]
    or
    \[
        \tau=\tau_2
        \quad\text{and}\quad
        \phi_2\in
        \left\{
        \frac{\pi_\beta}{2},
        \frac{3\pi_\beta}{2}
        \right\}.
    \]
\end{lemma}

\begin{proof}
    For \(w_0\neq0\), the claim follows from Theorem~\ref{Determinant}. We treat the boundary
    cases \(w_0=0\).

    First assume \(w_0=0\) and \((u_0,v_0)\neq(\pm1,0)\). In spherical coordinates this
    corresponds to \(\phi_2\to0\) or \(\phi_2\to\pi_\beta\), with
    \(\phi_1\notin\{0,\pi_\alpha\}\). We consider the limit as \(\phi_2\to0\); the case
    \(\phi_2\to\pi_\beta\) is identical. Let
    \[
        I(t):=\int_0^t x^{2\alpha}(s)\,ds,
        \qquad
        \varepsilon:=\frac{\delta_1}{y_0}.
    \]
    Then
    \[
        \omega_2=\varepsilon\sin_\beta(\phi_2),
        \qquad
        B=\phi_2+\varepsilon\sin_\beta(\phi_2)I(t).
    \]
    Using the Taylor expansions of \(\sin_\beta\) and \(\cos_\beta\) at \(0\), we have
    \[
        \sin_\beta(B)
        \sim
        \sin_\beta(\phi_2)\bigl(1+\varepsilon I(t)\bigr),
        \qquad
        \cos_\beta(B)\sim 1,
        \qquad
        \phi_2\to0.
    \]
    The \(\phi_2\)-factor in \eqref{det formula 3D} has the indeterminate form
    \(0/|\sin_\beta(\phi_2)|^{\beta+1}\). Differentiating numerator and denominator gives
    \[
        \begin{aligned}
             & \frac{d}{d\phi_2}
            \left[
                \cos_\beta(B)
                \left(\cot_\beta(\phi_2)\omega_2 I(t)+1\right)
                -
                \cot_\beta(\phi_2)\sin_\beta(B)
                \right]
            \\
             & \qquad\sim
            -\beta\sin_\beta^{2\beta-1}(\phi_2)
            \left[
                (1+\varepsilon I(t))^3
                +
                \varepsilon I(t)
                -
                (1+\varepsilon I(t))
                \right],
        \end{aligned}
    \]
    while
    \[
        \frac{d}{d\phi_2}
        |\sin_\beta(\phi_2)|^{\beta+1}
        \sim
        (\beta+1)|\sin_\beta(\phi_2)|^\beta .
    \]
    Thus the spherical determinant degenerates by a factor
    \(|\sin_\beta(\phi_2)|^{\beta-1}\). This is exactly the degeneracy of the spherical
    coordinate system. Indeed,
    \[
        \left|
        \frac{\partial(u_0,v_0,w_0)}
        {\partial(t,\phi_1,\phi_2)}
        \right|
        =
        \frac{
            t^2\alpha\beta
            \sin_\alpha^{2\alpha-1}(\phi_1)
            |\sin_\beta(\phi_2)|^{\beta-1}
        }{
            x_0^{2\alpha}|y_0|^\beta
        }.
    \]
    Therefore the Cartesian determinant
    \[
        D^{\operatorname{Cart}}\exp_{q_0}
        =
        \left|
        \frac{\partial(u_0,v_0,w_0)}
        {\partial(t,\phi_1,\phi_2)}
        \right|^{-1}
        D(t,\phi_1,\phi_2)
    \]
    has a finite limit as \(\phi_2\to0\). More precisely,
    \begin{equation}\label{e.cart-limit-w0}
        \begin{aligned}
            \lim_{\phi_2\to0}
            D^{\operatorname{Cart}}\exp_{q_0}
             & =
            C(\phi_1,t)
            \left[
                \cos_\alpha(A)
                \left(\cot_\alpha(\phi_1)\omega_1t+1\right)
                -
                \cot_\alpha(\phi_1)\sin_\alpha(A)
                \right]
            \\
             & \quad\times
            \left[
                (1+\varepsilon I(t))^3
                +
                \varepsilon I(t)
                -
                (1+\varepsilon I(t))
                \right],
        \end{aligned}
    \end{equation}
    where \(C(\phi_1,t)\neq0\) for \(t>0\) and
    \(\phi_1\notin\{0,\pi_\alpha\}\). The last bracket is
    \[
        (1+\varepsilon I)^3+\varepsilon I-(1+\varepsilon I)
        =
        \varepsilon I\bigl((1+\varepsilon I)^2+(1+\varepsilon I)+2\bigr).
    \]
    Since \(I(t)>0\) for \(t>0\), and \(\varepsilon\neq0\), this factor has no positive zero.
    Consequently, the only possible zeros of the limiting Cartesian determinant are the zeros
    of
    \[
        \cos_\alpha(A)
        \left(\cot_\alpha(\phi_1)\omega_1t+1\right)
        -
        \cot_\alpha(\phi_1)\sin_\alpha(A).
    \]
    By Proposition~\ref{tan function solution}, this factor does not vanish for
    \(0<t<\tau_1\). Since \(\tau_2=+\infty\) when \(w_0=0\), there are no conjugate times in
    \((0,\tau)\) for \(w_0=0\), \((u_0,v_0)\neq(\pm1,0)\).

    It remains to consider the straight-line geodesics \(\lambda_0=(\pm1,0,0)\). Write
    \(\eta=\pm1\). Then
    \[
        x(t)=x_0+\eta t,
        \qquad
        y(t)=y_0,
        \qquad
        z(t)=z_0.
    \]
    Use local Cartesian coordinates \((t,v_0,w_0)\) on the energy level near
    \((\eta,0,0)\), with \(u_0\) determined by
    \[
        u_0^2+x_0^{2\alpha}v_0^2+x_0^{2\alpha}y_0^{2\beta}w_0^2=1.
    \]
    At \(v_0=w_0=0\), one has
    \[
        \frac{\partial x}{\partial v_0}
        =
        \frac{\partial x}{\partial w_0}
        =
        0,
    \]
    \[
        \frac{\partial y}{\partial v_0}
        =
        \int_0^t (x_0+\eta s)^{2\alpha}\,ds,
        \qquad
        \frac{\partial y}{\partial w_0}=0,
    \]
    and
    \[
        \frac{\partial z}{\partial w_0}
        =
        y_0^{2\beta}
        \int_0^t (x_0+\eta s)^{2\alpha}\,ds,
        \qquad
        \frac{\partial z}{\partial v_0}=0.
    \]
    Since \(\partial_t x=\eta\), we obtain
    \begin{equation}\label{e.straight-line-det}
        \det
        \frac{\partial(x,y,z)}
        {\partial(t,v_0,w_0)}
        =
        \eta\,y_0^{2\beta}
        \left(
        \int_0^t (x_0+\eta s)^{2\alpha}\,ds
        \right)^2.
    \end{equation}
    This is nonzero for every \(t>0\). Hence straight-line geodesics have no conjugate
    points, which agrees with \(\tau_1=\tau_2=+\infty\).

    The equality statement follows from Theorem~\ref{Determinant} in the non-boundary cases
    and from the preceding limiting computation when \(w_0=0\).
\end{proof}

\section{Extended Hadamard technique for \texorpdfstring{$\mathbb{G}^3_{(\alpha,\beta)}$}{the three-dimensional Grushin space}}
Lemma~\ref{Con Times 3D} identifies the possible cut-conjugate points. We next determine,
on the slice \(H_{q_0}^{-1}(1/2)\cap\{u_0=0\}\), whether the smaller candidate time is
\(\tau_1\) or \(\tau_2\). This distinction controls the geometry of the two candidate
pieces of the cut locus.

\begin{proposition}\label{Tau1 bigger than Tau2}
    Let \(q_0=(x_0,y_0,z_0)\) be a Riemannian point in
    \(\mathbb G^3_{(\alpha,\beta)}\). Then
    \[
        \tau_1(\lambda_0)\leq \tau_2(\lambda_0)
        \qquad
        \text{for all }
        \lambda_0\in H_{q_0}^{-1}(1/2)\cap\{u_0=0\}
    \]
    if and only if
    \begin{equation}\label{separated condition}
        |y_0|
        \geq
        \frac{\pi_\alpha |x_0|^{\alpha+1}}
        {\pi_\beta(\alpha+1)}.
    \end{equation}
    If instead
    \[
        H_{q_0}^{-1}(1/2)\cap\{u_0=0\}\cap\{\tau_2\leq\tau_1\}
        \neq\emptyset,
    \]
    then the ellipse \(H_{q_0}^{-1}(1/2)\cap\{u_0=0\}\) splits into four arcs: two
    on which \(\tau_1\leq\tau_2\), centered at \(w_0=0\), and two on which
    \(\tau_2\leq\tau_1\), centered at the points of maximal \(|w_0|\).
\end{proposition}

\begin{proof}
    On \(H_{q_0}^{-1}(1/2)\), the condition \(u_0=0\) is
    \begin{equation}\label{e.u0-ellipse}
        v_0^2+y_0^{2\beta}w_0^2
        =
        |x_0|^{-2\alpha}.
    \end{equation}
    For \(u_0\neq\pm1\), the parameters in the geodesic formula satisfy
    \begin{equation}\label{parameter formulae}
        |\omega_1|
        =
        (v_0^2+y_0^{2\beta}w_0^2)^{1/(2\alpha)},
        \qquad
        A_1^{2\alpha}
        =
        (v_0^2+y_0^{2\beta}w_0^2)^{-1},
    \end{equation}
    and
    \[
        |\omega_2|
        =
        |w_0|^{1/\beta}
        (v_0^2+y_0^{2\beta}w_0^2)^{(\beta-1)/(2\beta)}.
    \]
    Thus
    \[
        \tau_2\leq\tau_1
    \]
    if and only if
    \[
        \pi_\beta
        \leq
        |\omega_2|\int_0^{\tau_1}x^{2\alpha}(s)\,ds.
    \]
    Using Proposition~\ref{Recursion Result} with \(\phi_1=\pi_\alpha/2\), this becomes
    \[
        \pi_\beta
        \leq
        \frac{|\omega_2|A_1^{2\alpha}}{|\omega_1|(\alpha+1)}\pi_\alpha.
    \]
    Equivalently,
    \begin{equation}\label{Inequality}
        w_0^2
        \geq
        \left(
        \frac{\pi_\beta(\alpha+1)}{\pi_\alpha}
        \right)^{2\beta}
        \bigl(v_0^2+y_0^{2\beta}w_0^2\bigr)^{\beta+1+\beta/\alpha}.
    \end{equation}
    Restricting \eqref{Inequality} to the ellipse \eqref{e.u0-ellipse}, we obtain
    \[
        w_0^2
        \geq
        \left(
        \frac{\pi_\beta(\alpha+1)}{\pi_\alpha}
        \right)^{2\beta}
        |x_0|^{-2(\alpha\beta+\alpha+\beta)}.
    \]
    Hence the set where \(\tau_2\leq\tau_1\) is the intersection of the ellipse
    \eqref{e.u0-ellipse} with the two horizontal regions
    \[
        |w_0|\geq
        \left(
        \frac{\pi_\beta(\alpha+1)}{\pi_\alpha}
        \right)^\beta
        |x_0|^{-(\alpha\beta+\alpha+\beta)}.
    \]
    Such an intersection is nonempty if and only if this lower bound is no larger than the
    maximal value of \(|w_0|\) on the ellipse, namely
    \[
        |w_0|_{\max}
        =
        \frac{1}{|x_0|^\alpha |y_0|^\beta}.
    \]
    Therefore \(\tau_2\leq\tau_1\) occurs somewhere on the slice \(u_0=0\) if and only if
    \[
        \frac{1}{|x_0|^\alpha |y_0|^\beta}
        \geq
        \left(
        \frac{\pi_\beta(\alpha+1)}{\pi_\alpha}
        \right)^\beta
        |x_0|^{-(\alpha\beta+\alpha+\beta)},
    \]
    which is equivalent to
    \[
        |y_0|
        \leq
        \frac{\pi_\alpha |x_0|^{\alpha+1}}
        {\pi_\beta(\alpha+1)}.
    \]
    This proves the first claim.

    The arc decomposition follows from the same description. On the ellipse
    \eqref{e.u0-ellipse}, the inequality \(\tau_2\leq\tau_1\) is equivalent to
    \(|w_0|\) lying above a fixed positive threshold. Thus, when the intersection is nonempty,
    there are two arcs with \(\tau_2\leq\tau_1\), centered at the points of maximal
    \(|w_0|\), and two complementary arcs with \(\tau_1\leq\tau_2\), centered at \(w_0=0\).
\end{proof}
\begin{figure}
    \centering
    \includegraphics[width=0.6\linewidth]{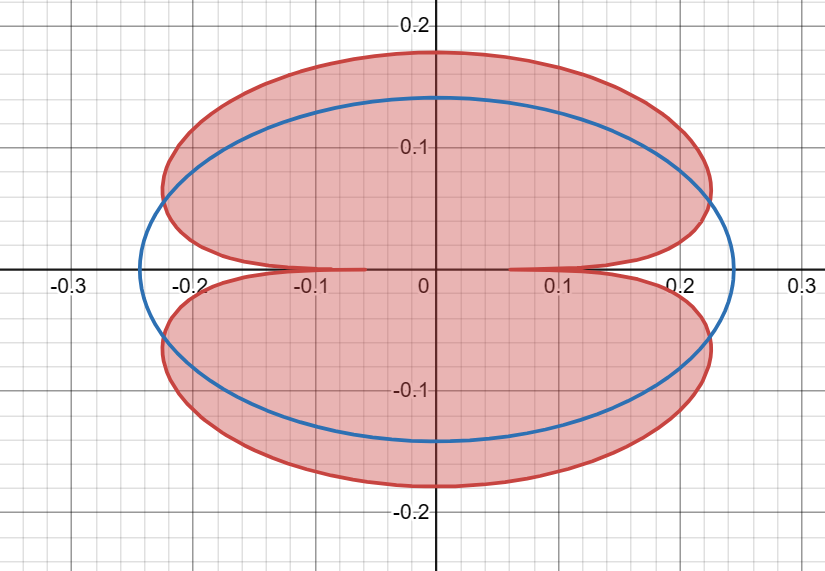}
    \caption{The red region in the $v_0w_0$-plane where $\{P(v_0,w_0)\geq 0\}$ or equivalently where \eqref{Inequality} holds, and the blue curve $H^{-1}_{q_0}(1/2)\cap\{u_0=0\}=\{y_0^{2\beta}w_0^2+v_0^2=\frac{1}{x^{2\alpha}}\}$ in the case when $q_0$ is an \emph{overlap point}.}
    \label{Inequality 1 and ellipse}
\end{figure}

\begin{definition}\label{point types}
    Let \(q_0=(x_0,y_0,z_0)\) be a Riemannian point in
    \(\mathbb G^3_{(\alpha,\beta)}\). We call \(q_0\) \emph{separated} if
    \begin{equation}\label{separated point definition}
        |y_0|
        \geq
        \frac{\pi_\alpha |x_0|^{\alpha+1}}
        {\pi_\beta(\alpha+1)}.
    \end{equation}
    We call \(q_0\) an \emph{overlap point} if
    \begin{equation}\label{overlap point definition}
        |y_0|
        \leq
        \frac{\pi_\alpha |x_0|^{\alpha+1}}
        {\pi_\beta(\alpha+1)}.
    \end{equation}
    If the corresponding inequality is strict, we say \(q_0\) is \emph{strictly separated}
    or \emph{strictly overlapping}, respectively.
\end{definition}
For separated points, the candidate endpoint loci corresponding to \(\tau_1\) and
\(\tau_2\) have disjoint boundaries inside
\[
    \{x=-x_0\}\cup\{y=-y_0\}.
\]
For overlap points, these boundaries meet. The boundary points are precisely the
candidate cut-conjugate points identified in Lemma~\ref{Con Times 3D}. This distinction
will be used in the proof of Theorem~\ref{main theorem}.

\begin{figure}
    \centering
    \includegraphics[scale=.5]{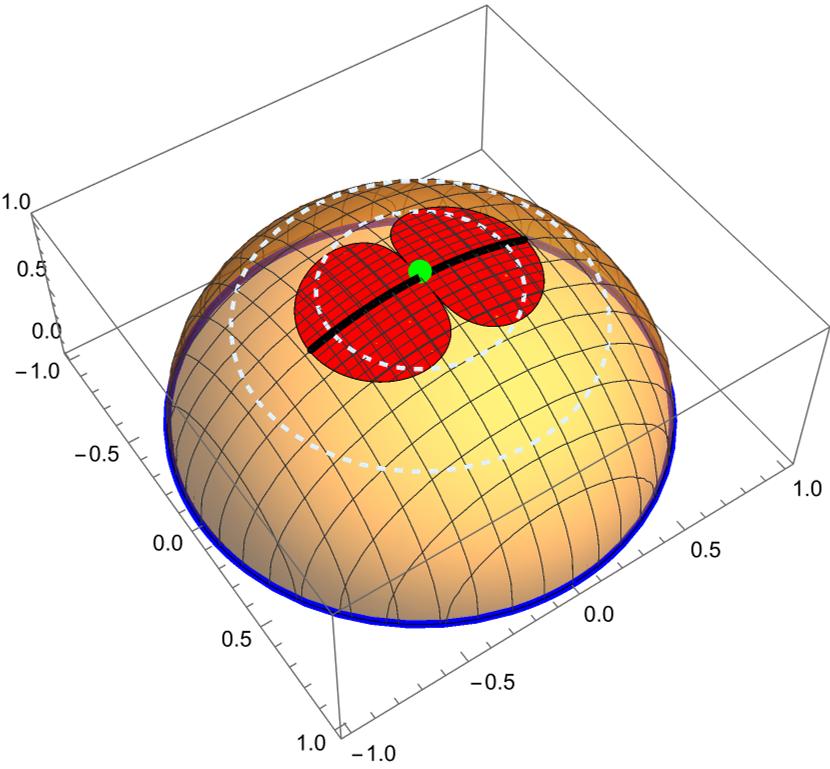}
    \caption{The set $H^{-1}_{q_0}(1/2)$ for $q_0$ a \emph{separated Riemannian point}, oriented in with the $u_0$ coordinate facing upward, and $(1,0,0)$ the \textcolor{green} {thick green dot}. The \textcolor{red}{red region} represents $\{\tau_2\leq\tau_1\}.$ The thick black curve represents $\Lambda_2=\{\lambda_0\in H^{-1}_{q_0}: \tau_2(\lambda_0)\leq \tau_1(\lambda_0),v_0=0\}.$ Observe that the thick \textcolor{blue}{blue curve} representing $H^{-1}_{q_0}(1/2)\cap\{u_0=0\}$, consists entirely of points in $\{\tau_1\leq\tau_2\}$}
    \label{fig:Tau1 and Tau2 Regions}
\end{figure}

We next describe the boundary of the candidate endpoint locus corresponding to
\(\tau_2\). By Lemma~\ref{Con Times 3D}, \(\tau_2\) is conjugate precisely when
\[
    \phi_2\in
    \left\{
    \frac{\pi_\beta}{2},
    \frac{3\pi_\beta}{2}
    \right\}
    \qquad\text{and}\qquad
    \tau_2\leq\tau_1.
\]
Set
\[
    \Lambda_2
    :=
    \left\{
    \lambda_0\in H_{q_0}^{-1}(1/2):
    \phi_2\in
    \left\{
    \frac{\pi_\beta}{2},
    \frac{3\pi_\beta}{2}
    \right\},
    \ \tau_2\leq\tau_1
    \right\}.
\]
For \(\lambda_0\in\Lambda_2\), the endpoint
\(\exp_{q_0}(\tau_2(\lambda_0)\lambda_0)\) lies in the plane \(y=-y_0\). Moreover,
\[
    \partial_{\phi_1}z
    \left(
    \tau_2,\phi_1,\frac{\pi_\beta}{2}
    \right)
    =
    0,
\]
and similarly for \(\phi_2=3\pi_\beta/2\). Thus, along each component of \(\Lambda_2\),
the \(z\)-coordinate is constant. For \(\phi_2=\pi_\beta/2\), using the recursion formula
for \(z\),
\[
    \begin{aligned}
        z(\tau_2)
         & =
        z_0+
        \frac{w_0A_2^{2\beta}}{\omega_2(\beta+1)}
        \left[
            \eta_\beta
            \left(
            \omega_2\int_0^{\tau_2}x^{2\alpha}(s)\,ds+\frac{\pi_\beta}{2}
            \right)
            -
            \eta_\beta\left(\frac{\pi_\beta}{2}\right)
            \right]
        \\
         & =
        z_0+
        \frac{w_0A_2^{2\beta}}{\omega_2(\beta+1)}
        \left[
            \eta_\beta\left(\frac{3\pi_\beta}{2}\right)
            -
            \eta_\beta\left(\frac{\pi_\beta}{2}\right)
            \right]
        \\
         & =
        z_0+
        \frac{\pi_\beta |y_0|^{\beta+1}}{\beta+1}.
    \end{aligned}
\]
The case \(\phi_2=3\pi_\beta/2\) gives the negative sign. Hence the candidate
\(\tau_2\)-boundary is contained in the two horizontal lines
\[
    \left\{
    (x,-y_0,z_0\pm
    \tfrac{\pi_\beta |y_0|^{\beta+1}}{\beta+1})
    :
    x\in\mathbb R
    \right\}.
\]
In the proof of Theorem~\ref{main theorem}, we show that these lines are exactly the
\(\tau_2\)-boundary:
\[
    \left\{
    \exp_{q_0}(\tau_2(\lambda_0)\lambda_0):
    \lambda_0\in\Lambda_2
    \right\}
    =
    \left\{
    (x,-y_0,z_0\pm
    \tfrac{\pi_\beta |y_0|^{\beta+1}}{\beta+1})
    :
    x\in\mathbb R
    \right\}.
\]

We will also use the following elementary consequence of the inverse function theorem. See \cite{LeeBook2003SmoothManifold}.

\begin{lemma}\label{Crit Points}
    Let \(M,N\) be smooth manifolds of the same dimension, \(U\subset M\) open, and
    \(h:U\to N\) a \(C^1\) map. If \(q\in U\) and \(h(q)\in\partial h(U)\), then $\ker d_qh\neq\{0\}$.
\end{lemma}
\begin{proof}
    If \(d_qh\) were invertible, the inverse function theorem would imply that \(h(U)\)
    contains a neighborhood of \(h(q)\), contradicting \(h(q)\in\partial h(U)\).
\end{proof}

\begin{figure}
    \centering
    \includegraphics[scale=.85]{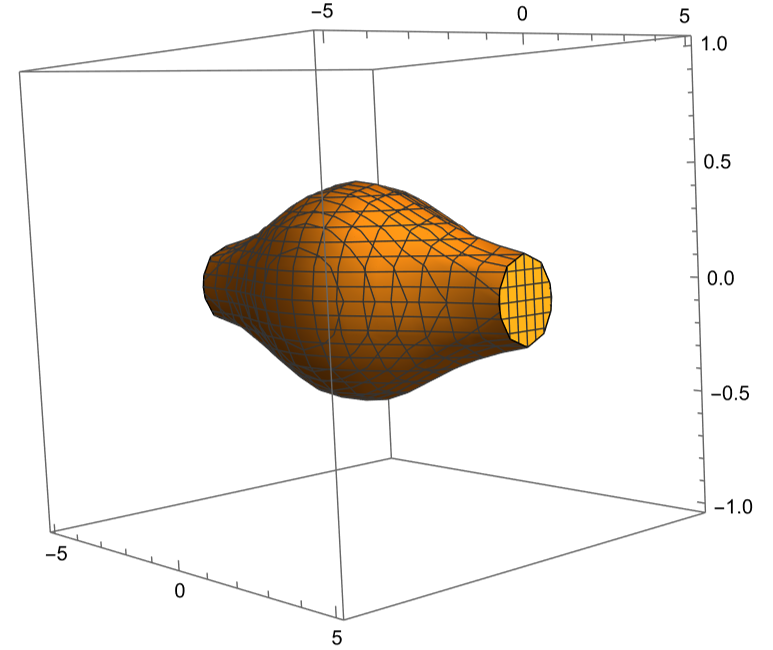}
    \caption{A superset \(N_1\) of the cotangent injectivity domain \(N\), defined by
        \(|\omega_1|<\pi_\alpha\). It is unbounded only in the \(u_0\)-direction.}
    \label{fig:Cotangent Injectivity}
\end{figure}

\begin{figure}
    \centering
    \includegraphics[scale=.5]{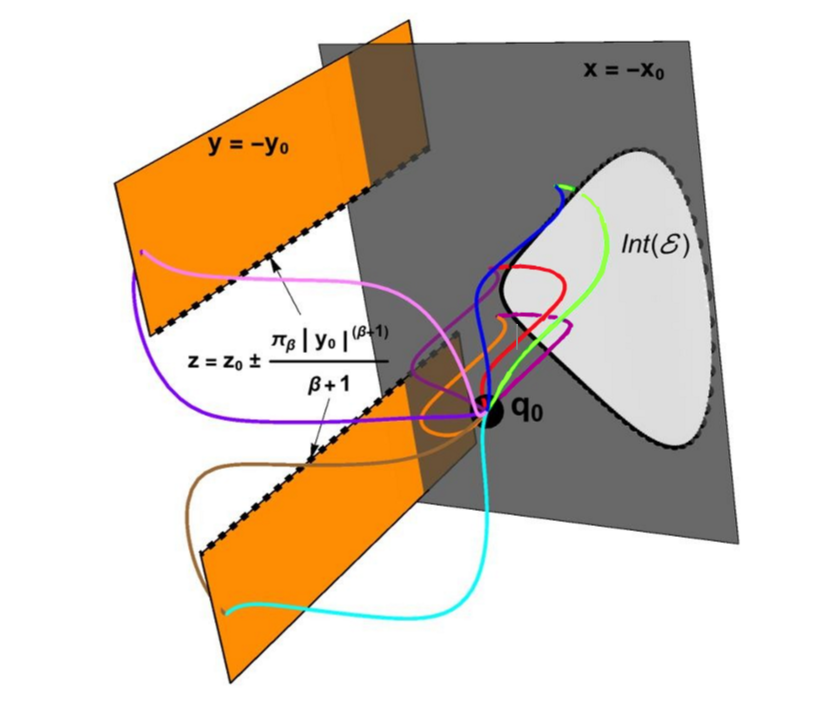}
    \caption{Candidate cut locus for a separated Riemannian point
        \(q_0=(1,1,z_0)\) in \(\mathbb G^3_{(\alpha,\beta)}\), with
        \(\alpha=\beta=2\). The dotted horizontal lines
        \(z=z_0\pm \pi_\beta |y_0|^{\beta+1}/(\beta+1)\) and the curve
        \(\mathcal E\) contain the candidate cut-conjugate points.}
    \label{fig:Cut Locus}
\end{figure}

\begin{theorem}\label{main theorem}
    Let \(\alpha,\beta\geq1\), and let
    \(q_0=(x_0,y_0,z_0)\in\mathbb G^3_{(\alpha,\beta)}\) be Riemannian. If
    \(\gamma(t)=\exp_{q_0}(t\lambda_0)\) is arclength-parametrized, then
    \[
        t_{\operatorname{cut}}(\gamma)=\tau(\lambda_0).
    \]
    Moreover,
    \[
        \operatorname{Cut}(q_0)
        =
        \operatorname{Cut}_x^*(q_0)
        \cup
        \operatorname{Cut}_y^*(q_0),
    \]
    where the $\tau_2$ locus of endpoints is given by
    \[
        \operatorname{Cut}_y^*(q_0)
        :=
        \left\{
        (x,-y_0,z):
        |z-z_0|
        \geq
        \frac{\pi_\beta |y_0|^{\beta+1}}{\beta+1}
        \right\},
    \]
    while the $\tau_1$ locus of endpoints is
    \[
        \operatorname{Cut}_x^*(q_0)
        :=
        \{x=-x_0\}\setminus \operatorname{int}(\mathcal E).
    \]
    and we define
    \[
        \mathcal E
        :=
        \left\{
        \exp_{q_0}(\tau_1(\lambda_0)\lambda_0):
        \lambda_0\in \Lambda_1
        \right\},
    \]
    \[
        \Lambda_1
        :=
        \left\{
        \lambda_0\in H_{q_0}^{-1}(1/2):
        u_0=0,\ \tau_1(\lambda_0)\leq \tau_2(\lambda_0)
        \right\},
    \]
    and \(\operatorname{int}(\mathcal E)\) denotes the bounded component of
    \(\{x=-x_0\}\setminus\mathcal E\).
\end{theorem}

\begin{proof}
    Set
    \[
        \operatorname{Cut}^*(q_0)
        :=
        \operatorname{Cut}_x^*(q_0)
        \cup
        \operatorname{Cut}_y^*(q_0),
    \]
    and
    \[
        N
        :=
        \left\{
        t\lambda_0:
        \lambda_0\in H_{q_0}^{-1}(1/2),
        \ 0\leq t<\tau(\lambda_0)
        \right\}.
    \]
    We verify the hypotheses of Theorem~\ref{Extended Hadamard Technique}. We will verify each of i)–iii) in the assumptions of Theorem \ref{Extended Hadamard Technique} separately. We proceed by verifying ii), then iii), and finally i).\smallskip

    \textbf{Proof of ii)} We must show that $\exp_{q_0}\rvert_N$ is a proper, local diffeomorphism. First, \(\exp_{q_0}\) has no critical points on $N$ by
    Lemma~\ref{Con Times 3D}. It remains to check properness. Since
    \[
        N\subset N_1:=\{|\omega_1|<\pi_\alpha\},
    \]
    the variables \(v_0,w_0\) are bounded on \(N_1\), while \(u_0\) is the only unbounded
    direction. If \(|u_0|\to\infty\) with \(|\omega_1|<\pi_\alpha\), then
    \[
        u_0=A_1\omega_1\cos_\alpha(\phi_1)
    \]
    forces \(|A_1|\to\infty\). Hence the \(x\)-coordinate of
    \(\exp_{q_0}(u_0,v_0,w_0)\), which is controlled by $\lvert A_1\rvert$, tends to infinity. Thus
    \(\exp_{q_0}|_N\) is proper.

    \textbf{Proof of iii)} We must show that $\mathbb{G}^3_{(\alpha,\beta)}\setminus \operatorname{Cut}^*(q_0)$ is simply connected. We proceed by studying the topology of \(\mathcal E\).

    \begin{claim}\label{simply connected claim}
        The set \(\mathcal E\subset\{x=-x_0\}\) is a simple closed curve
    \end{claim}
    In the separated case, Proposition~\ref{Tau1 bigger than Tau2} gives
    \[
        \Lambda_1
        =
        \left\{
        \lambda_0\left(\frac{\pi_\alpha}{2},\phi_2\right):
        \phi_2\in[0,2\pi_\beta)
        \right\},
    \]
    so \(\mathcal E\) is a closed parametrized curve.

    In the overlap case, let \(\phi_2^*\) be the first parameter for which
    \(\tau_1=\tau_2\) on the slice \(u_0=0\). Then \(\mathcal E\) is parametrized by
    \[
        [0,\phi_2^*]
        \cup
        [\pi_\beta-\phi_2^*,\pi_\beta+\phi_2^*]
        \cup
        [2\pi_\beta-\phi_2^*,2\pi_\beta],
    \]
    with the endpoints glued by the symmetry
    \(\phi_2\mapsto\pi_\beta-\phi_2\) and
    \(\phi_2\mapsto2\pi_\beta-\phi_2\).

    It remains only to rule out self-intersections. Along the \(\tau_1\)-endpoint map,
    \[
        \partial_{\phi_2}y(\tau_1,\pi_\alpha/2,\phi_2)
        =
        \frac{y_0}{\sin_\beta(\phi_2)}
        \left[
            \cos_\beta(B)
            \left(
            \cot_\beta(\phi_2)\omega_2
            \int_0^{\tau_1}x^{2\alpha}(s)\,ds
            +1
            \right)
            -
            \cot_\beta(\phi_2)\sin_\beta(B)
            \right].
    \]
    If this derivative vanished away from
    \(\phi_2=\pi_\beta/2,3\pi_\beta/2\), then Proposition~\ref{tan function solution}
    would imply
    \[
        \left|
        \omega_2\int_0^{\tau_1}x^{2\alpha}(s)\,ds
        \right|
        >
        \pi_\beta,
    \]
    contradicting \(\tau_1\leq\tau_2\) on \(\Lambda_1\). Hence \(y\) is monotone on
    each relevant arc. The symmetries of \(\sin_\beta\) and \(\cos_\beta\), together with
    \(z\geq z_0\) on the upper half and \(z\leq z_0\) on the lower half, rule out
    intersections between distinct arcs. Thus \(\mathcal E\) is simple and closed, which completes the proof of Claim \ref{simply connected claim}.

    By the Jordan curve theorem, \(\mathcal E\) bounds a topological disk in the plane
    \(\{x=-x_0\}\). Therefore \(\operatorname{Cut}_x^*(q_0)\) is that plane with the
    bounded disk removed, while \(\operatorname{Cut}_y^*(q_0)\) consists of two closed
    half-planes in \(\{y=-y_0\}\). After a translation and a homeomorphism preserving the
    coordinate planes, the complement
    \[
        \mathbb G^3_{(\alpha,\beta)}\setminus\operatorname{Cut}^*(q_0)
    \]
    is homeomorphic to \(\mathbb R^3\) with two properly embedded half-planes and one
    properly embedded punctured plane removed. This complement is simply connected; for
    example, one cuts along the plane \(\{y=-y_0\}\) and applies van Kampen to the two
    resulting simply connected regions, whose intersection is connected and simply
    connected. This completes the proof of iii).

    \textbf{Proof of i)}. It is required to show the equality
    \begin{align}\label{main set equality}
        \operatorname{exp}_{q_0}(N)=\mathbb{G}^3_{(\alpha,\beta)}\setminus \operatorname{Cut}^*(q_0).
    \end{align}
    We now proceed by a boundary identification argument of the image sets under two component endpoint maps.

    We first attempt to show
    \begin{align}\label{claim 2}
        F_1\bigl(\{\tau_1\leq\tau_2\}\setminus\{(\pm1,0,0)\}\bigr)
        =
        \operatorname{Cut}_x^*(q_0),
    \end{align}
    where
    \[
        F_1(\lambda_0):=\exp_{q_0}(\tau_1(\lambda_0)\lambda_0).
    \]
    Indeed, \(F_1\) takes values in the plane \(\{x=-x_0\}\). It is proper: as
    \(\lambda_0\to(\pm1,0,0)\), the \((y,z)\)-coordinates of \(F_1(\lambda_0)\) diverge.
    On \(\{\tau_1<\tau_2\}\setminus\{u_0=0\}\), a determinant computation similar to what was done in Section 6 shows that
    \(F_1\) is a local diffeomorphism. The set \(u_0=0\) is fixed by the symmetry
    \(u_0\mapsto-u_0\), and \(F_1\) is invariant under this symmetry. Hence
    \(\Lambda_1\setminus\{\tau_1=\tau_2\}\) is a fold. By the standard fold normal form,
    its image is a boundary arc of \(F_1(\{\tau_1\leq\tau_2\})\). Therefore
    \[
        F_1(\Lambda_1)=\mathcal E
        \subset
        \partial F_1(\{\tau_1\leq\tau_2\}).
    \]
    The remaining boundary comes from \(\{\tau_1=\tau_2\}\), whose image lies on
    \[
        \{x=-x_0,\ y=-y_0\}
    \]
    and joins the endpoints of the corresponding arcs of \(\mathcal E\). These are exactly
    the cut-conjugate boundary points already identified in Lemma~\ref{Con Times 3D}.
    Thus the boundary of the image is \(\mathcal E\). Since \(F_1\) is proper, its image is
    closed; since it is unbounded, it must be the exterior of \(\mathcal E\) in the plane
    \(\{x=-x_0\}\). Hence
    \[
        F_1(\{\tau_1\leq\tau_2\})
        =
        \{x=-x_0\}\setminus\operatorname{int}(\mathcal E)
        =
        \operatorname{Cut}_x^*(q_0).
    \]
    which proves \eqref{claim 2}. We prove an analogous equation for the endpoint map under $\tau_2$:
    \begin{align}\label{claim 3}F_2\bigl(\{\tau_2\leq\tau_1\}\bigr)
        =
        \operatorname{Cut}_y^*(q_0)\end{align}
    where
    \[
        F_2(\lambda_0):=\exp_{q_0}(\tau_2(\lambda_0)\lambda_0).
    \]
    The proof is analogous but simpler. The map \(F_2\) takes values in the plane
    \(\{y=-y_0\}\). It is invariant under the symmetry \(v_0\mapsto -v_0\), and the critical
    set in \(\{\tau_2<\tau_1\}\) is precisely \(\{v_0=0\}\). Thus the image of
    \(\Lambda_2\) is the boundary of the image. By the computation preceding
    Lemma~\ref{Crit Points},
    \[
        F_2(\Lambda_2)
        =
        \left\{
        (x,-y_0,z_0\pm
        \frac{\pi_\beta |y_0|^{\beta+1}}{\beta+1})
        :
        x\in\mathbb R
        \right\}.
    \]
    Since \(F_2\) is also proper and its image is unbounded in the \(z\)-direction, the image is
    the union of the two closed half-planes
    \[
        \operatorname{Cut}_y^*(q_0)
        =
        \left\{
        (x,-y_0,z):
        |z-z_0|
        \geq
        \frac{\pi_\beta |y_0|^{\beta+1}}{\beta+1}
        \right\}.
    \]
    \begin{figure}
        \centering
        \includegraphics[scale=.65]{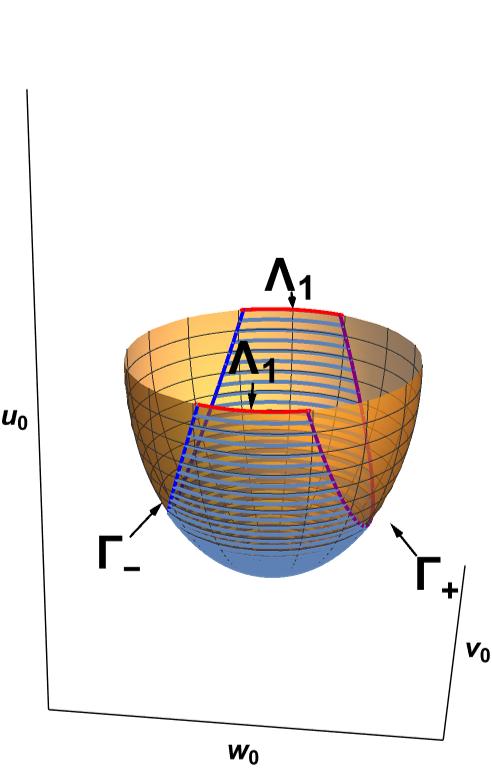}
        \caption{Illustration of the blue shaded region $\mathcal{T}\subset H^{-1}_{q_0}(1/2)$ in the overlap case, where $\mathcal{T}=\{t^*\leq \tau\}$,  with boundary given by $\Lambda_1$ and the curves $\Gamma_{-},\Gamma_{+}$.}
        \label{Mathcal T Region}
    \end{figure}
    We are ready now to prove \eqref{main set equality}.

    For the inclusion \(\supseteq\) in \eqref{main set equality}, let
    \(q\notin\operatorname{Cut}^*(q_0)\). Since the space is complete and has no strictly
    abnormal minimizers from a Riemannian point \(q_0\), \(q\) is joined to \(q_0\) by a normal minimizing
    geodesic \(\exp_{q_0}(t\lambda_0)\). By Theorem~\ref{upper bound cut time},
    \(t\leq\tau(\lambda_0)\). Equality would put \(q\) in
    \(\operatorname{Cut}^*(q_0)\), by Claims \eqref{claim 2} and \eqref{claim 3}. Hence \(t<\tau(\lambda_0)\), and
    \(q\in\exp_{q_0}(N)\).

    For the inclusion \(\subseteq\), let
    \[
        q=\exp_{q_0}(t\lambda_0),
        \qquad
        0<t<\tau(\lambda_0).
    \]
    We must show \(q\notin\operatorname{Cut}^*(q_0)\). The only possible issue is when
    \(q\) lies in one of the planes \(x=-x_0\) or \(y=-y_0\).

    Suppose first \(q\in\{x=-x_0\}\). Then \(q\) is reached at the first hitting time
    \[
        t^*
        =
        \frac{2(\pi_\alpha-\phi_1)}{|\omega_1|},
        \qquad
        \frac{\pi_\alpha}{2}<\phi_1<\pi_\alpha.
    \]
    Let
    \[
        G_1(\lambda_0)
        :=
        \exp_{q_0}(t^*(\lambda_0)\lambda_0).
    \]
    On the region
    \[
        \mathcal T
        :=
        \{\lambda_0:t^*(\lambda_0)\leq\tau(\lambda_0)\},
    \]
    the boundary consists of \(\Lambda_1\), and in the overlap case, the two arcs
    \(\Gamma_\pm\) where \(t^*=\tau_2\). The boundary identity
    \[
        |\sin_\beta(\phi_2)|K(\phi_1)
        =
        \frac{\pi_\beta |y_0|(\alpha+1)}{|x_0|^{\alpha+1}},
        \qquad
        K(\phi_1)
        =
        \frac{2(\pi_\alpha-\eta_\alpha(\phi_1))}
        {\sin_\alpha^{\alpha+1}(\phi_1)},
    \]
    shows that
    \[
        G_1(\Gamma_\pm)
        \subset
        \operatorname{int}(\mathcal E)
        \cap
        \operatorname{Cut}_y^*(q_0).
    \]
    Moreover, \(G_1(\Lambda_1)=\mathcal E\), and \(G_1\) is a local diffeomorphism in the
    interior of \(\mathcal T\). Therefore
    \[
        G_1(\{t^*<\tau\})
        \subset
        \operatorname{int}(\mathcal E)\setminus \operatorname{Cut}_y^*(q_0).
    \]
    Thus \(q\notin\operatorname{Cut}_x^*(q_0)\cup\operatorname{Cut}_y^*(q_0)\).

    Now suppose \(q\in\{y=-y_0\}\). Then \(q\) is reached at the first hitting time
    \(t^{**}\), determined by
    \[
        \omega_2\int_0^{t^{**}}x^{2\alpha}(s)\,ds
        =
        2\pi_\beta-2\phi_2,
        \qquad
        \phi_2\in\left[\frac{\pi_\beta}{2},\frac{3\pi_\beta}{2}\right].
    \]
    At this time,
    \[
        \begin{aligned}
            |z(t^{**})-z_0|
             & =
            \frac{|y_0|^{\beta+1}}
            {|\sin_\beta(\phi_2)|^{\beta+1}}
            \int_{\phi_2}^{2\pi_\beta-\phi_2}
            \sin_\beta^{2\beta}(s)\,ds.
        \end{aligned}
    \]
    This function is zero at \(\phi_2=\pi_\beta\), increases toward
    \(\pi_\beta/2\), and by symmetry decreases on the other side. Its maximum is attained
    at
    \[
        \phi_2=\frac{\pi_\beta}{2},
        \qquad
        \phi_2=\frac{3\pi_\beta}{2},
    \]
    with value
    \[
        \frac{\pi_\beta |y_0|^{\beta+1}}{\beta+1}.
    \]
    Since \(t^{**}<\tau\), equality cannot occur. Hence
    \[
        |z(t^{**})-z_0|
        <
        \frac{\pi_\beta |y_0|^{\beta+1}}{\beta+1},
    \]
    so \(q\notin \operatorname{Cut}_y^*(q_0)\).

    It remains to exclude \(q\in\operatorname{Cut}_x^*(q_0)\). If the point is overlapping,
    then \(t^{**}\leq\tau_1\) throughout the relevant sector, and this is immediate. In the
    separated case, one repeats the preceding first-hit argument with the roles of \(x\) and
    \(y\) interchanged: on
    \[
        \mathcal T'
        :=
        \{\lambda_0:t^{**}(\lambda_0)\leq\tau(\lambda_0)\},
    \]
    the endpoint map
    \[
        G_2(\lambda_0)
        :=
        \exp_{q_0}(t^{**}(\lambda_0)\lambda_0)
    \]
    has boundary contained in the boundary of \(\operatorname{int}(\mathcal E)\), and its
    interior maps away from \(\{x=-x_0\}\). Hence \(q\notin\operatorname{Cut}_x^*(q_0)\).
    Therefore \(q\notin\operatorname{Cut}^*(q_0)\). Both inclusions in the set equality \eqref{main set equality} are then proved.

    The Extended Hadamard Technique now applies: \(\exp_{q_0}|_N\) is a proper local
    diffeomorphism onto a simply connected target, hence a diffeomorphism. Therefore no
    geodesic loses optimality before \(\tau\). Together with
    Theorem~\ref{upper bound cut time}, this gives
    \[
        t_{\operatorname{cut}}(\gamma)=\tau(\gamma).
    \]
    The cut locus is consequently \(\operatorname{Cut}^*(q_0)\).
\end{proof}
\begin{figure}
    \centering
    \includegraphics[scale=.5]{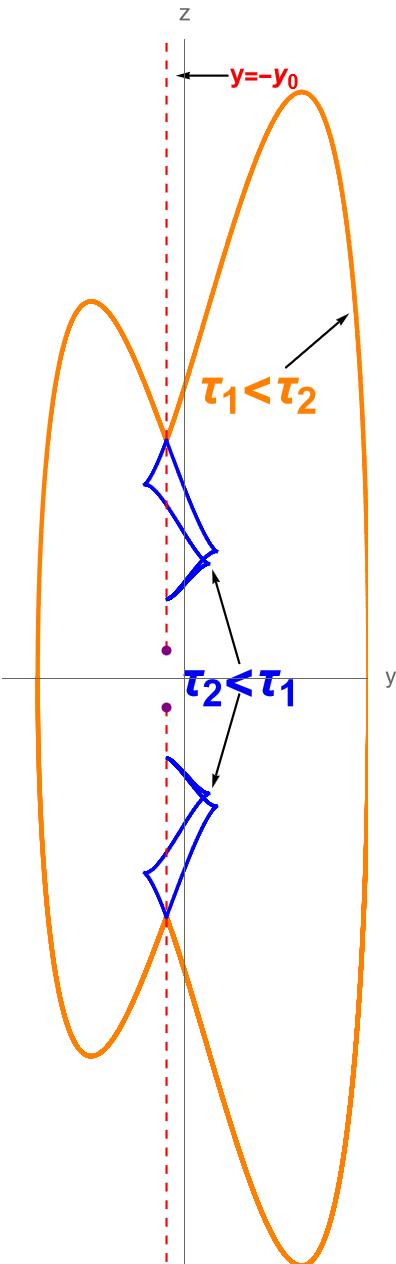}
    \caption{For the overlapping point $q_0=(2.4,1,0)$ and $(\alpha,\beta)=(1,1)$, this is the image $\exp_{q_0}(\tau_1(\lambda_0)\lambda_0)$ in the $yz$-plane for $\lambda_0\in H^{-1}_{q_0}(1/2)\cap\{u_0=0\}$. The \textcolor{orange}{orange portion of the curve} is by definition $\mathcal{E}$, which we observe is simple and closed. The \textcolor{blue}{blue portion with $\tau_2<\tau_1$} has some self intersections.}
    \label{fig:enter-label}
\end{figure}
\section{Optimal synthesis at singular points by density}\label{Singular Points}
We now extend the optimal synthesis from Riemannian initial points to singular ones in
\(\mathbb G^3_{(\alpha,\beta)}\). The only additional issue is that the fibers
\(H_{q_0}^{-1}(1/2)\) are no longer smooth compact hypersurfaces when \(q_0\) is singular.
We use a density argument together with upper semicontinuity of the cut time.

We recall the following standard fact.

\begin{theorem}[{\cite[Prop. 8.76]{AgrachevBarilariBoscainBook2020}}]
    \label{cut time cont}
    The cut time function
    \[
        t_{\operatorname{cut}}:H^{-1}(1/2)\to[0,+\infty]
    \]
    is upper semicontinuous.
\end{theorem}

We first extend the candidate cut time to singular fibers.

\begin{lemma}\label{Singular Opt Synth}
    The function
    \[
        \tau=\min\{\tau_1,\tau_2\}
    \]
    extends continuously on \(H^{-1}(1/2)\) wherever it is finite. Its infinite values occur
    precisely on covectors whose normal geodesic is the horizontal line
    \[
        t\mapsto (x_0\pm t,y_0,z_0).
    \]
\end{lemma}

\begin{proof}
    Let \(q_0=(x_0,y_0,z_0)\). The Hamiltonian is
    \[
        H(q_0,u_0,v_0,w_0)
        =
        \frac12
        \left(
        u_0^2+x_0^{2\alpha}v_0^2+x_0^{2\alpha}y_0^{2\beta}w_0^2
        \right).
    \]
    At singular points the spherical coordinates used at Riemannian points may degenerate,
    but the frequencies entering \(\tau_1\) and \(\tau_2\) have continuous Cartesian
    expressions:
    \begin{equation}\label{e.singular-frequencies}
        |\omega_1|
        =
        \bigl(v_0^2+y_0^{2\beta}w_0^2\bigr)^{1/(2\alpha)},
        \qquad
        |\omega_2|
        =
        |w_0|^{1/\beta}
        \bigl(v_0^2+y_0^{2\beta}w_0^2\bigr)^{(\beta-1)/(2\beta)}.
    \end{equation}
    These formulae agree with the Riemannian ones and are obtained directly from the
    conserved quantities
    \[
        R_2^2=v_0^2+y_0^{2\beta}w_0^2,
        \qquad
        R_3^2=w_0^2.
    \]

    If \(R_2=0\), then \(x(t)=x_0\pm t\), while \(y(t)\equiv y_0\) and
    \(z(t)\equiv z_0\) whenever the corresponding transverse dynamics vanish. These are the
    horizontal line geodesics, and we set
    \[
        \tau_1=\tau_2=+\infty.
    \]

    If \(R_2>0\), then
    \[
        \tau_1=\frac{\pi_\alpha}{|\omega_1|}.
    \]
    The \(x\)-coordinate is the continuous solution of
    \[
        x''=-\alpha x^{2\alpha-1}R_2^2,
        \qquad
        x(0)=x_0,\quad x'(0)=u_0,
    \]
    and therefore depends continuously on \((q_0,\lambda_0)\). Hence
    \[
        \tau_2
        =
        \inf\left\{
        t>0:
        |\omega_2|\int_0^t x^{2\alpha}(s)\,ds=\pi_\beta
        \right\}
    \]
    is continuous wherever it is finite. Thus
    \[
        \tau=\min\{\tau_1,\tau_2\}
    \]
    is continuous on its finite locus.
\end{proof}

\begin{theorem}\label{full optimal synth}
    In \(\mathbb G^3_{(\alpha,\beta)}\),
    \[
        t_{\operatorname{cut}}=\tau
    \]
    on \(H^{-1}(1/2)\), with \(\tau\) defined as above.
\end{theorem}

\begin{proof}
    Let
    \[
        \Omega
        :=
        \bigcup_{q_0\notin\mathcal S}
        \{q_0\}\times H_{q_0}^{-1}(1/2)
        \subset H^{-1}(1/2),
    \]
    where \(\mathcal S=\{xy=0\}\) is the singular set. Then \(\Omega\) is dense in
    \(H^{-1}(1/2)\). By Theorem~\ref{main theorem},
    \[
        t_{\operatorname{cut}}=\tau
        \qquad
        \text{on }\Omega.
    \]

    The symmetry argument used in Theorem~\ref{upper bound cut time} applies also on singular
    fibers, using the singular geodesic formulas from Lemma~\ref{Singular Opt Synth}; hence
    \begin{equation}\label{e.global-upper-bound}
        t_{\operatorname{cut}}\leq \tau
        \qquad
        \text{on }H^{-1}(1/2).
    \end{equation}

    Now fix \(\lambda\in H^{-1}(1/2)\) with \(\tau(\lambda)<+\infty\), and choose
    \(\lambda_k\in\Omega\) with \(\lambda_k\to\lambda\). By Lemma~\ref{Singular Opt Synth},
    \[
        \tau(\lambda_k)\to\tau(\lambda).
    \]
    Since \(t_{\operatorname{cut}}(\lambda_k)=\tau(\lambda_k)\), upper semicontinuity gives
    \[
        \tau(\lambda)
        =
        \lim_{k\to\infty}t_{\operatorname{cut}}(\lambda_k)
        \leq
        t_{\operatorname{cut}}(\lambda).
    \]
    Together with \eqref{e.global-upper-bound}, this yields
    \[
        t_{\operatorname{cut}}(\lambda)=\tau(\lambda).
    \]

    If \(\tau(\lambda)=+\infty\), then \(\lambda\) generates a horizontal line geodesic.
    Such a curve is globally minimizing, so
    \[
        t_{\operatorname{cut}}(\lambda)=+\infty=\tau(\lambda).
    \]
    Thus \(t_{\operatorname{cut}}=\tau\) on all of \(H^{-1}(1/2)\).
\end{proof}

\begin{remark}
    The same density strategy would apply in higher-dimensional triangular Grushin spaces
    once Conjecture~\ref{conjecture} is proved. The essential input is the absence of
    strictly abnormal minimizers via Theorem \ref{abnormality}: every unit-speed length minimizer admits a normal lift in
    \(H^{-1}(1/2)\), even when it is also abnormal.
\end{remark}
We conclude with the description of the cut locus at singular initial points.
\begin{corollary}\label{cor:singular-cut-locus}
    Let \(q_0=(x_0,y_0,z_0)\in \mathbb G^3_{(\alpha,\beta)}\) be singular. Then
    \[
        \operatorname{Cut}(q_0)=
        \begin{cases}
            (\{x=0\}\cup(\{y=0\})\setminus\{(x,0,z_0)\})\,
             & x_0=0,\ y_0=0,     \\[1mm]
            \{x=0\}\cup
            \Bigl(
            \{y=-y_0\}
            \cap
            \bigl\{|z-z_0|\geq\dfrac{\pi_\beta |y_0|^{\beta+1}}{\beta+1}
            \bigr\}
            \Bigr),
             & x_0=0,\ y_0\neq 0, \\[2mm]
            \Bigl(\{x=-x_0\}\setminus \operatorname{int}(\mathcal G)\Bigr)\cup (\{y=0\}\setminus\{(x,0,z_0)\}),
             & x_0\neq 0,\ y_0=0,
        \end{cases}
    \]
    where
    \[
        \mathcal G
        :=
        \left\{
        \exp_{q_0}(\tau_1(\lambda_0)\lambda_0):
        u_0=0,\ \tau_1(\lambda_0)\leq \tau_2(\lambda_0)
        \right\}
    \]
    is a simple closed curve in the plane \(\{x=-x_0\}\). Equivalently, in Cartesian
    coordinates on \(H^{-1}_{q_0}(1/2)\),
    \[
        \mathcal G
        =
        \left\{
        \exp_{q_0}(\tau_1(\lambda_0)\lambda_0):
        u_0=0,\
        v_0=\pm \frac{1}{|x_0|^\alpha},\
        |w_0|
        \leq
        \left(
        \frac{\pi_\beta}{\pi_\alpha}(\alpha+1)
        \right)^\beta
        \frac{1}{|x_0|^{\alpha\beta+\alpha+\beta}}
        \right\}.
    \]
\end{corollary}

\begin{proof}
    By Theorem~\ref{full optimal synth}, every normal geodesic satisfies
    \(t_{\mathrm{cut}}=\tau=\min\{\tau_1,\tau_2\}\), with the convention that
    \(t_{\mathrm{cut}}=+\infty\) when \(\tau=+\infty\). Hence the cut locus is obtained by taking
    only those endpoints for which the cut time is finite:
    \[
        \operatorname{Cut}(q_0)
        =
        F_1\bigl(\{\tau_1\leq \tau_2,\ \tau_1<+\infty\}\bigr)
        \cup
        F_2\bigl(\{\tau_2\leq \tau_1,\ \tau_2<+\infty\}\bigr),
    \]
    where
    \[
        F_1(\lambda_0)=\exp_{q_0}(\tau_1(\lambda_0)\lambda_0),
        \qquad
        F_2(\lambda_0)=\exp_{q_0}(\tau_2(\lambda_0)\lambda_0).
    \]
    The covectors for which \(\tau=+\infty\) generate precisely the horizontal straight-line
    geodesics
    \[
        t\longmapsto (x_0+\eta t,y_0,z_0),\qquad \eta=\pm1.
    \]
    These geodesics are globally minimizing, and therefore their images do not contribute to
    \(\operatorname{Cut}(q_0)\).

    We now identify the two finite endpoint images case by case.

    First suppose \(x_0=0\) and \(y_0=0\). The \(\tau_1\)-endpoints lie in the plane
    \(\{x=0\}\), while the \(\tau_2\)-endpoints lie in the plane \(\{y=0\}\). There are no
    cut-conjugate points in either finite family. Consequently, away from the straight-line
    covectors, the maps \(F_1\) and \(F_2\) are proper local diffeomorphisms onto their images.
    The same open-and-closed argument used in the proof of Theorem \ref{main theorem} gives
    \[
        F_1\bigl(\{\tau_1\leq \tau_2,\ \tau_1<+\infty\}\bigr)=\{x=0\},
    \]
    and
    \[
        F_2\bigl(\{\tau_2\leq \tau_1,\ \tau_2<+\infty\}\bigr)
        =
        \{y=0\}\setminus\{(x,0,z_0):x\in\mathbb R\}.
    \]
    The excluded line is exactly the image of the globally minimizing straight-line geodesics.
    Thus
    \[
        \operatorname{Cut}(q_0)
        =
        (\{x=0\}\cup\{y=0\})
        \setminus
        \{(x,0,z_0):x\in\mathbb R\}.
    \]

    Next suppose \(x_0=0\) and \(y_0\neq0\). The finite \(\tau_1\)-endpoints fill the plane
    \(\{x=0\}\). The finite \(\tau_2\)-endpoints lie in the plane \(\{y=-y_0\}\). The condition
    \(t=\tau_2\) gives
    \[
        |z-z_0|\geq \frac{\pi_\beta |y_0|^{\beta+1}}{\beta+1}.
    \]
    Conversely, every point in this set is reached by a finite \(\tau_2\)-endpoint. Hence
    \[
        F_2\bigl(\{\tau_2\leq \tau_1,\ \tau_2<+\infty\}\bigr)
        =
        \{y=-y_0\}
        \cap
        \left\{
        |z-z_0|
        \geq
        \frac{\pi_\beta |y_0|^{\beta+1}}{\beta+1}
        \right\}.
    \]
    This gives the second formula.

    Finally suppose \(x_0\neq0\) and \(y_0=0\). The finite \(\tau_2\)-endpoints lie in
    \(\{y=0\}\). As in the first case, the straight-line covectors have \(\tau=+\infty\) and
    their image is the line
    \[
        \{(x,0,z_0):x\in\mathbb R\}.
    \]
    Removing this globally minimizing line, the open-and-closed argument gives
    \[
        F_2\bigl(\{\tau_2\leq \tau_1,\ \tau_2<+\infty\}\bigr)
        =
        \{y=0\}\setminus\{(x,0,z_0):x\in\mathbb R\}.
    \]

    For the \(\tau_1\)-family, the endpoints lie in the plane \(\{x=-x_0\}\). The
    cut-conjugate covectors are precisely those satisfying
    \[
        u_0=0,\qquad \tau_1(\lambda_0)\leq \tau_2(\lambda_0),
    \]
    or, equivalently,
    \[
        u_0=0,\qquad
        v_0=\pm \frac{1}{|x_0|^\alpha},\qquad
        |w_0|
        \leq
        \left(
        \frac{\pi_\beta}{\pi_\alpha}(\alpha+1)
        \right)^\beta
        \frac{1}{|x_0|^{\alpha\beta+\alpha+\beta}}.
    \]
    Their image under \(F_1\) is the simple closed curve \(\mathcal G\). The same monotonicity
    argument as in the Riemannian case shows that \(\mathcal G\) is simple and closed. Away
    from this critical set, \(F_1\) is a proper local diffeomorphism, and the symmetry
    \(u_0\mapsto -u_0\) shows that the critical set is a fold. Hence
    \[
        \mathcal G\subset \partial\operatorname{Im}(F_1).
    \]
    Since \(\operatorname{Im}(F_1)\subset\{x=-x_0\}\) is the unbounded component determined by
    \(\mathcal G\), we obtain
    \[
        F_1\bigl(\{\tau_1\leq\tau_2,\ \tau_1<+\infty\}\bigr)
        =
        \{x=-x_0\}\setminus\operatorname{int}(\mathcal G).
    \]
    Therefore
    \[
        \operatorname{Cut}(q_0)
        =
        \Bigl(\{x=-x_0\}\setminus\operatorname{int}(\mathcal G)\Bigr)
        \cup
        \Bigl(\{y=0\}\setminus\{(x,0,z_0):x\in\mathbb R\}\Bigr),
    \]
    which proves the corollary.
\end{proof}
\begin{figure}
    \centering
    \includegraphics[width=0.5\linewidth]{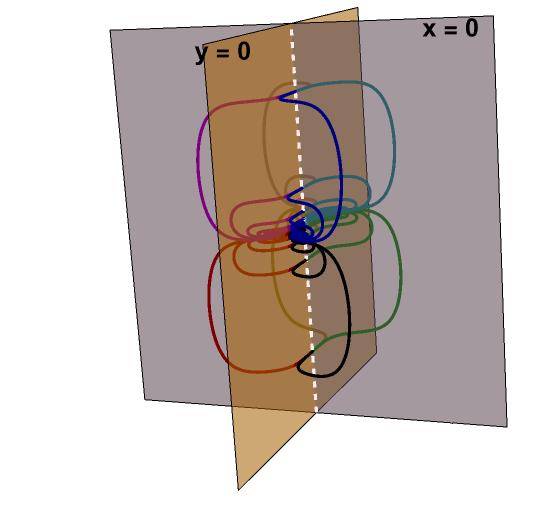}
    \caption{$\operatorname{Cut}(q_0)$ for a singular point $q_0$ with $x_0=y_0=0$ consisting of the plane $\{y=0\}$ together with $\{x=0\}$. The family of geoedesics displayed are chosen to intersect on the dotted $z$-axis, where $\tau_1=\tau_2$. There are no cut-conjugate points.}
    \label{fig:placeholder}
\end{figure}

\begin{figure}

    \centering
    \includegraphics[width=.5\textwidth]{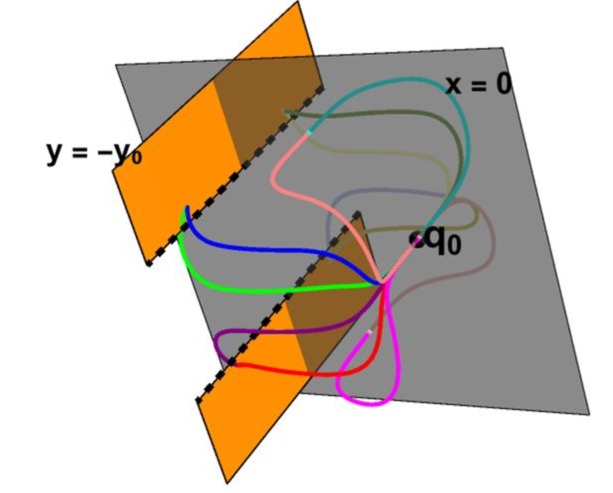}
    \caption{$\operatorname{Cut}(q_0)$ for a singular point $q_0$ with $x_0=0,y_0\neq 0$ consisting of the plane $\{x=0\}$ together with the two strips $\{\lvert z-z_0\rvert\geq \frac{\pi_\beta\lvert y_0\rvert^{\beta+1}}{\beta+1}\}\cap\{y=-y_0\}$. The cut-conjugate points lie on the dotted lines.}
    \label{fig:fig8}

\end{figure}
\begin{figure}

    \centering
    \includegraphics[width=.5\textwidth]{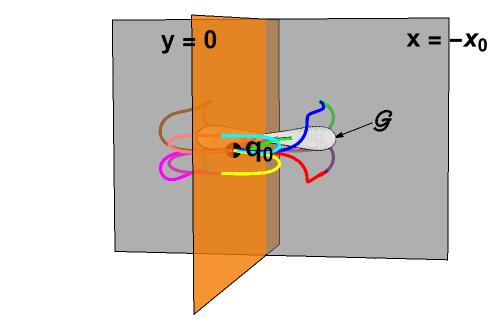}
    \caption{$\operatorname{Cut}(q_0)$ for a singular point $q_0$ with $x_0\neq 0,y_0=0$ consisting of the plane $\{y=0\}$ together with $\{x=-x_0\}\setminus\operatorname{int}(\mathcal{G})$. The cut-conjugate points lie on $\mathcal{G}$.}
    \label{fig:fig9}

\end{figure}
\newpage

\bibliographystyle{plain}
\bibliography{Grushin}

\end{document}